\documentclass[12 pt]{article}
\usepackage[hmargin=1in,vmargin=1in]{geometry}
\usepackage[latin1]{inputenc}
\usepackage{subfigure}
\usepackage{color}
\usepackage{graphicx}
\usepackage{amsfonts}
\usepackage{amsmath}
\usepackage{graphics}
\usepackage{amssymb}
\usepackage{latexsym}
\usepackage{amscd}
\usepackage{color}
\usepackage{pict2e}
\usepackage[all]{xy}
\usepackage{bm}
\usepackage{enumerate}
\usepackage{todonotes}
\usepackage{mathrsfs}
\usepackage{float}
\usepackage{amssymb,amsthm,amsmath,amsfonts,latexsym,tikz,hyperref,color,enumitem}
\newtheorem{thm}{Theorem}[section]
\newtheorem{prop}[thm]{Proposition}
\newtheorem{lem}[thm]{Lemma}
\newtheorem{cor}[thm]{Corollary}

\newtheorem{question}[thm]{Question}
\theoremstyle{definition}

\newtheorem{definition}[thm]{Definition}
\newtheorem{algorithm}[thm]{Algorithm}

\newcommand{\NC}{{NC}}
\newcommand{\NBB}{\mathcal{NBB}}

\newcommand{\NCNBC}{\mathcal{NCNBC}}

\newcommand{\cover}{\lessdot}
\renewcommand{\c}{\mathbf{c}}
\DeclareMathOperator{\ncisf}{ncisf}
\DeclareMathOperator{\ncnbc}{ncnbc}

\DeclareMathOperator{\nbc}{nbc}

\DeclareMathOperator{\dist}{dist}

\title {The  Noncrossing Bond Poset of a Graph}

\author{C. Matthew Farmer\\
\small Department of Mathematics and Statistics\\[-0.8ex]
\small The University of North Carolina at Greensboro\\[-0.8ex]
\small Greensboro, NC, 27412, U.S.A.\\
\small \tt cmfarmer@uncg.edu\\
\and
Joshua  Hallam\\
\small Department of Mathematics\\[-0.8ex]
\small Loyola Marymount University\\[-0.8ex]
\small Los Angeles, CA 90045, U.S.A.\\
\small \tt joshua.hallam@lmu.edu\\
\and
Clifford Smyth\thanks{Supported by Simons Collaboration Grant 360486.}\\
\small Department of Mathematics and Statistics\\[-0.8ex]
\small The University of North Carolina at Greensboro\\[-0.8ex]
\small Greensboro, NC, 27412, U.S.A.\\
\small \tt  cdsmyth@uncg.edu
}

\begin{document}

\maketitle

\begin{abstract}
    
The partition lattice and noncrossing partition lattice   are well studied objects in combinatorics.  Given a graph $G$ on vertex set $\{1,2,\dots, n\}$, its bond lattice, $L_G$, is the subposet of the partition lattice formed by restricting to the partitions whose blocks induce connected subgraphs of $G$.  In this article, we introduce a natural noncrossing analogue of the bond lattice, the noncrossing bond poset, $NC_G$, obtained by  restricting to the noncrossing partitions of $L_G$.

Both the noncrossing partition lattice and the bond lattice have many nice combinatorial properties.  We show that, for several families of graphs, the noncrossing bond poset also exhibits these properties. We present simple necessary and sufficient conditions on the graph to ensure the noncrossing bond poset is a lattice.  Additionally, for several families of graphs, we give combinatorial descriptions of the M\"obius function and characteristic polynomial of the noncrossing bond poset. These descriptions are in terms of a noncrossing analogue of non-broken circuit (NBC) sets of the graphs and can be thought of as a noncrossing version of Whitney's NBC theorem for the chromatic polynomial. We also consider the shellability and supersolvability of the noncrossing bond poset, providing sufficient conditions for both.   We end with some open problems. 
    

\end{abstract}

\setcounter{tocdepth}{1}
\tableofcontents

\section{Introduction}\label{introSec}

The noncrossing partition lattice on $[n] = \{1,2, \dots, n\}$, denoted by $\NC_n$ in this paper, was introduced in 1972 by Kreweras~\cite{k:nc} and has received considerable attention from the combinatorics community since then.   It enjoys  many nice order-theoretic properties.  For example, it is graded, EL-shellable~\cite{b:scmpos},   self-dual (and hence rank-symmetric)~\cite{k:nc}, $k$-Sperner for all $k$~\cite{su:slnc},  and supersolvable~\cite{h:depos}. The Catalan numbers appear in several contexts of the noncrossing partition lattice, being both the number of elements of the lattice as well as the M\"obius value of the maximum element of the lattice. The maximal chains of $\NC_n$ are equinumerous with spanning trees of the complete graph $K_n$, and thus also  with parking functions. In~\cite{s:pfnc}, Stanley showed a beautiful connection between parking functions and $\NC_n$ by exhibiting a natural edge labeling indexed by parking functions.  Consequently, Stanley showed that the flag (quasi)symmetric function of $\NC_n$ is (up to a simple automorphism) the parking function symmetric function introduced by Haiman~\cite{h:cqrdi}.  In addition to their significance in combinatorics, noncrossing partitions have surprising connections to  other mathematics including low-dimensional topology, geometric group theory, mathematical biology,  and noncommutative probability.  We refer the reader to McCammond's~\cite{m:ncsl} and   Simion's~\cite{s:np}  survey articles for more information on the many properties of noncrossing partitions both within and outside of combinatorics.

The noncrossing partition lattice has been generalized in several ways.   In~\cite{r:nccfg}, Reiner introduced a type-B and type-D noncrossing partition lattice and showed that they have many of the nice properties that $\NC_n$ (the type-A version) has.  Reiner then asked if one could define a noncrossing partition lattice for any Coxeter group.  The work of Bessis~\cite{b:dbm}, Brady~\cite{b:posym}, and Brady and Watt~\cite{bw:kPi} showed that not only is this possible, but that these lattices play an important role in constructing monoids and $K(\pi,1)$ spaces for the Artin groups associated with finite Coxeter groups.  One recovers $\NC_n$ by taking the symmetric group as the Coxeter group.  This family of noncrossing partition lattices associated to Coxeter groups retains many of the nice order-theoretic properties that $\NC_n$ possesses.  For example, they are graded, lattices (see~\cite{bw:ncfrrg} for a uniform proof), shellable (see~\cite{abw:snc} for a uniform proof), and self-dual~\cite{b:dbm}.  

Other methods of generalizing the noncrossing partition lattice have been studied and shown to be fruitful.  For example, in~\cite{e:cenc} Edelman introduced the  $k$-divisible noncrossing partition lattice, the subposet of $\NC_n$ where each partition has all of its block sizes divisible by $k$.  Later, Armstrong~\cite{a:gnc} introduced and studied the $k$-divisible noncrossing partition lattice for each finite Coxeter group. More recently, motivated by the connection between the noncrossing partition lattice and parking functions,  Bruce, Dougherty, Hlavacek, Kudo, and Nicolas~\cite{bdhkn:dpf} introduced a subposet of the noncrossing partition lattice obtained by removing chains that corresponded to parking functions with certain restrictions.  To solve a conjecture put forward in that article, M\"uhle~\cite{m:tpncps} defined two new subposets of the noncrossing partition lattice obtained by removing partitions which do not contain certain blocks.  M\"uhle~\cite{m:tpncps}  showed that these new posets are graded, shellable, and supersolvable.

In this article, we introduce a new generalization of the noncrossing partition lattice which is based on the structure of finite graphs. This generalization can be thought of as the intersection of the noncrossing partition lattice and a bond lattice.    Given a graph with vertex set $[n]$, its bond lattice is a subposet of the partition lattice obtained by restricting it to the set of partitions such that for each block $B$ in the partition, the induced subgraph of $G$ with vertex set $B$ is connected. Note that the bond lattice of the complete graph is the partition lattice since any induced subgraph of the complete graph is connected.  The bond lattice carries important combinatorial  information about the graph.   For example, it encodes exactly the same information as the cycle matroid associated to the graph. In fact, the bond lattice is (isomorphic to) the lattice of flats of this cycle matroid.   Moreover, its characteristic polynomial is (essentially) the chromatic polynomial of the graph and its chromatic symmetric function can be computed from the lattice as well.   Since the partition lattice is the bond lattice of the complete graph, one can consider the noncrossing partition lattice as a noncrossing version of a bond lattice.  It is this idea that is the starting point of our work.   The \emph{noncrossing bond poset} of a graph is the intersection of the noncrossing partition lattice and its bond lattice, i.e.~the poset obtained from its bond lattice by removing any partition that is crossing.

While many of the generalizations of the noncrossing partition lattice discussed above  exhibit the nice properties of $\NC_n$, the situation for the noncrossing bond poset is a bit more nuanced.  In general, many of these properties do not hold for the noncrossing bond poset of generic graphs.  We note that this might be expected as the structure of graphs can vary widely.  Nevertheless, we are able to identify families of graphs for which some of these nice properties still hold.  We present simple necessary and sufficient conditions on the graph for the noncrossing bond poset to be a lattice (see Theorem~\ref{crossingClosedThm}).    The  M\"obius function and the characteristic polynomial of the bond lattice of a graph can be interpreted combinatorially in terms of non-broken circuit (NBC) sets via Whitney's NBC theorem.  We show that, for several families of graphs, similar interpretations hold for the noncrossing bond poset in terms of what we call noncrossing NBC sets. We obtain our noncrossing version of Whitney's NBC theorem in two different ways, one by using the theory of non-bounded below (NBB) sets introduced by Blass and Sagan~\cite{bs:mfl} (see Theorem~\ref{NCNBCThm}) and the other by using the minimum EL-labeling for geometric lattices introduced by Bj\"orner~\cite{b:scmpos} (see Theorem~\ref{minLabelNC}).  Both approaches are necessary as there are some graphs which can be handled by only one of these methods.  Moreover, the EL-labeling approach allows us to provide shellability results for the noncrossing bond poset of several families of graphs.     Additionally, we show that the noncrossing bond poset of perfectly labeled graphs, which arise from chordal graphs, admit $S_n$ EL-labelings (see Proposition~\ref{minMaxPerfectLabel}).  Using a result of McNamara~\cite{m:els}, we get that when these are lattices they are supersolvable lattices.  These results on perfectly labeled graphs parallel the results for chordal graphs in realm of  bond lattices.

We also give two algorithms for non-crossing bond posets in the appendix.  Algorithm \ref{crossingclosedalg} determines if the noncrossing bond poset of a graph is a lattice.  Algorithm \ref{uppercrossingclosedalg} determines if a graph belongs to a family of graphs for which the noncrossing NBC interpretation of the M\"obius function and characteristic polynomial hold.   Our algorithms both run in time polynomial in $n$, the number of vertices of the graph.  This is of interest because brute-force algorithms that do not take advantage of the theory we develop can, in the worst case, take time super-exponential in $n$.

The rest of the paper is organized as follows.  In Section~\ref{structureSec} we discuss the basic structure of the noncrossing bond poset.  After this, we consider the M\"obius function and characteristic polynomial in Section~\ref{MobSec}.   Edge labelings are then studied in Section~\ref{edgeLabSec}.  Next, we look at the properties of several families of graphs in Section~\ref{famOfGraphsSec}.  Section~\ref{openProblemsSec} contains several open problems.  We finish with the appendix on algorithms we previously mentioned.

Some of these results appeared in an earlier version of this paper published in the Proceedings of Formal Power Series and Algebraic Combinatorics 2019~\cite{fh:ncp}.

\section{The Structure of the Noncrossing Bond Poset}\label{structureSec}

\begin{figure}
\begin{center}
\begin{tikzpicture}[scale=.7]
\coordinate (v1) at (0,1);
\coordinate (v2) at (.87,.5);
\coordinate (v3) at (.87,-.5);
\coordinate (v4) at (0,-1);
\coordinate (v5) at (-.87,-.5);
\coordinate (v6) at (-.87,.5);

\draw[thick] (v6)--(v1) -- (v2) -- (v3)--(v5);
\draw[thick] (v5)--(v6);
\draw[thick] (v1) -- (v4);

\foreach \v in {v1,v2,v3,v4, v5,v6} \fill(\v) circle (.1);
\draw(0,1.4) node{\footnotesize 1};
\draw (.87+.4,.5+.4) node{\footnotesize 2};
\draw(0,-1.4) node{\footnotesize 4};
\draw (.87+.4,-.5-.4) node{\footnotesize 3};
\draw (-.87-.4,-.5-.4) node{\footnotesize 5};
\draw (-.87-.4,.5+.4) node{\footnotesize 6};
\draw(0,-2.2) node {\large $G$};

\begin{scope}[shift = {(4, 0)}]

\coordinate (v1) at (0,1);
\coordinate (v2) at (.87,.5);
\coordinate (v3) at (.87,-.5);
\coordinate (v4) at (0,-1);
\coordinate (v5) at (-.87,-.5);
\coordinate (v6) at (-.87,.5);

\draw[thick] (v1) -- (v2) -- (v3) --(v5);
\draw[thick] (v1) -- (v4);

\foreach \v in {v1,v2,v3,v4, v5,v6} \fill(\v) circle (.1);
\draw(0,1.4) node{\footnotesize 1};
\draw (.87+.4,.5+.4) node{\footnotesize 2};
\draw(0,-1.4) node{\footnotesize 4};
\draw (.87+.4,-.5-.4) node{\footnotesize 3};
\draw (-.87-.4,-.5-.4) node{\footnotesize 5};
\draw (-.87-.4,.5+.4) node{\footnotesize 6};
\draw(0,-2.2) node {\large $H$};

\end{scope}

\begin{scope}[shift = {(8, 0)}]

\coordinate (v1) at (0,1);
\coordinate (v2) at (.87,.5);
\coordinate (v3) at (.87,-.5);
\coordinate (v4) at (0,-1);
\coordinate (v5) at (-.87,-.5);
\coordinate (v6) at (-.87,.5);

\draw[thick] (v3) --(v5)--(v6);
\draw[thick] (v1) -- (v4);
\draw[thick] (v3) -- (v5);
\draw[thick] (v1) -- (v6)--(v5);

\foreach \v in {v1,v2,v3,v4, v5,v6} \fill(\v) circle (.1);
\draw(0,1.4) node{\footnotesize 1};
\draw (.87+.4,.5+.4) node{\footnotesize 2};
\draw(0,-1.4) node{\footnotesize 4};
\draw (.87+.4,-.5-.4) node{\footnotesize 3};
\draw (-.87-.4,-.5-.4) node{\footnotesize 5};
\draw (-.87-.4,.5+.4) node{\footnotesize 6};
\draw(0,-2.2) node {\large $K$};

\end{scope}

\begin{scope}[shift = {(12, 0)}]

\coordinate (v1) at (0,1);
\coordinate (v2) at (.87,.5);
\coordinate (v3) at (.87,-.5);
\coordinate (v4) at (0,-1);
\coordinate (v5) at (-.87,-.5);
\coordinate (v6) at (-.87,.5);

\draw[thick] (v1) -- (v4);
\draw[thick] (v3) -- (v5);

\foreach \v in {v1,v2,v3,v4, v5,v6} \fill(\v) circle (.1);
\draw(0,1.4) node{\footnotesize 1};
\draw (.87+.4,.5+.4) node{\footnotesize 2};
\draw(0,-1.4) node{\footnotesize 4};
\draw (.87+.4,-.5-.4) node{\footnotesize 3};
\draw (-.87-.4,-.5-.4) node{\footnotesize 5};
\draw (-.87-.4,.5+.4) node{\footnotesize 6};
\draw(0,-2.2) node {\large $H \cap K$};

\end{scope}

\begin{scope}[shift = {(16, 0)}]

\coordinate (v1) at (0,1);
\coordinate (v2) at (.87,.5);
\coordinate (v3) at (.87,-.5);
\coordinate (v4) at (0,-1);
\coordinate (v5) at (-.87,-.5);
\coordinate (v6) at (-.87,.5);

\draw[thick] (v1) -- (v2) -- (v3)--(v5)--(v6);
\draw[thick] (v1) -- (v4);

\foreach \v in {v1,v2,v3,v4, v5,v6} \fill(\v) circle (.1);
\draw(0,1.4) node{\footnotesize 1};
\draw (.87+.4,.5+.4) node{\footnotesize 2};
\draw(0,-1.4) node{\footnotesize 4};
\draw (.87+.4,-.5-.4) node{\footnotesize 3};
\draw (-.87-.4,-.5-.4) node{\footnotesize 5};
\draw (-.87-.4,.5+.4) node{\footnotesize 6};
\draw(0,-2.2) node {\large $L$};

\end{scope}

\end{tikzpicture}
\end{center}
\caption{A graph and several subgraphs \label{graphFig} }
\end{figure}
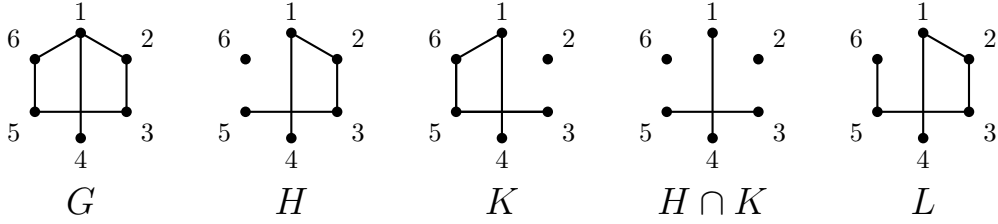

We assume the reader is familiar with the basic concepts of graph theory (see~\cite[Graph Theory Appendix]{s:ec1} for any undefined terms) as well as basic concepts related to posets (see~\cite[Chapter 3]{s:ec1} for background and notation).   Let $G$ be a graph. For the remainder of this paper, unless otherwise noted, we will assume that the vertex set of $G$ is $[n]$. We will use the notation $V(G)$ for the vertex set of $G$ and $E(G)$ for the edge set of $G$.   When we write out edges, we will write them in the form $ij$ where $i<j$.   Moreover, we will always draw our graphs so that the vertices lie on a circle with vertex $1$ at the top and the remaining vertices appearing in clockwise order around the circle. Edges will always be drawn so that they are  the line segments between their endpoints. We will refer to this as the \emph{graphical representation of $G$}.   We say that two edges of $G$ \emph{cross} if their respective line segments intersect in the graphical representation, i.e.~$a_1 a_2$ and $b_1 b_2$ cross if and only if $a_1 < b_1 < a_2 < b_2$ or $b_1 < a_1 < b_2 < a_2$. See Figure~\ref{graphFig} for several  examples of  graphical representations.

A subgraph $H$ of a graph $G$ is called \emph{spanning} if $V(H) = V(G)$. Note that when considering a spanning subgraph $H$ of $G$, it is enough to just know $E(H)$.  Because of this, we will often make no distinction between subsets of $E(G)$ and spanning subgraphs of $G$.    A subgraph $H$ of a graph $G$ is called \emph{induced} if whenever $u$ and $v$ are vertices in $H$ and $uv\in E(G)$, then $uv\in E(H)$. Given a subset of the vertices, $S$, let $G[S]$ denote the induced subgraph of $G$ with vertex set $S$.  Similarly, if $E$ is a set of edges, we will use $G[E]$ to denote the induced subgraph on the vertices which are endpoints of the edges in $E$.  We say a spanning subgraph of  $G$ is a \emph{bond} if every connected component of the subgraph is induced.  As an example, consider the graph $G$ in Figure~\ref{graphFig}. The subgraphs $H$, $K$, and $H \cap K$ are bonds, but $L$ is not since it is missing the edge $16$.  

To each bond $H$, one can associate a set partition, $\pi(H)$, so that $i$ and $j$ are in the same block of $\pi(H)$ if and only if $i$ and $j$ are in the same connected component of $H$.  For example,  for the bond $H$ in Figure~\ref{graphFig}, we have that $\pi(H)=12345/6$. Similarly, for each partition $\pi = B_1/B_2/ \cdots/B_k$, we associate the corresponding spanning subgraph $G[\pi]$ whose edge set is the disjoint union of the edges in $G[B_1], \dots, G[B_k]$.  We note that $G[B_1/B_2/ \cdots/B_k]$ is not necessarily a bond  since $G[B_i]$ need not be connected.

A partition $\pi = B_1/B_2/\cdots/ B_k$  is called \emph{crossing} if there exists $i \neq j$, $a,c \in B_i$ and $b,d \in B_j$ such that $a<b<c<d$.  For example, the partition  $1248/ 56/37$ is crossing since   we can pick $2, 4\in 1248$ and $3,7 \in 37$.  A partition is \emph{noncrossing} if it is not crossing.  We say the bond $H$ is  \emph{crossing} (resp.~\emph{noncrossing}) if $\pi(H)$ is crossing (resp.~noncrossing).   It is not hard to verify the following proposition. 

\begin{prop}\label{crossingBondProp}
A bond $H$ is crossing if and only if there exist two distinct connected components $H_1$ and $H_2$ of $H$ and edges $e_1 \in H_1$ and $e_2 \in H_2$ such that $e_1$ and $e_2$ cross.
\end{prop}

Note that a noncrossing bond $H$ can contain crossing edges, as long as every pair of edges that cross belong to the same connected component of $H$.     For example, the bond $H$ in Figure~\ref{graphFig} is noncrossing since it corresponds to $12345/6$, but it has crossing edges, namely $14$ and $35$.

\begin{definition}
Let $G$ be a graph.  The \emph{bond lattice} of $G$, denoted by $L_G$, is the collection of bonds of $G$ ordered by inclusion. The \emph{noncrossing bond poset}, denoted by $\NC_G$, is the collection of noncrossing  bonds of $G$ ordered by inclusion. 
\end{definition}

\begin{figure}
\begin{center}
\begin{tikzpicture}[scale=.7]
\coordinate (v1) at (0,1);
\coordinate (v2) at (1,0);
\coordinate (v3) at (0,-1);
\coordinate (v4) at (-1,0);
\draw[thick] (v3)--(v1);
\draw[thick] (v2) -- (v4);
\foreach \v in {v1,v2,v3,v4} \fill(\v) circle (.1);
\draw(0,1.4) node{\footnotesize 1};
\draw(0,-1.4) node{\footnotesize 3};
\draw(1.4,0) node{\footnotesize 2};
\draw(-1.4,0) node{\footnotesize 4};
\draw(0,-3.2) node {\large $G$};

\begin{scope}[shift = {(8,-2)}]

\node (1/2/3/4) at (0,0)  {$1/2/3/4$};
\node (13/2/4) at (-2,2)  {$13/2/4$};
\node (1/24/3) at (2,2)  {$1/24/3$};
\node (13/24) at (0,4) {$13/24$};

\draw[thick] (1/2/3/4)--(13/2/4);
\draw[thick] (1/2/3/4)--(1/24/3);

\draw[thick] (13/2/4)--(13/24);
\draw[thick] (1/24/3)--(13/24);

\draw(0,-1.2) node {\large $L_G$};

\end{scope}

\begin{scope}[shift = {(16,-2)}]

\node (1/2/3/4) at (0,0)  {$1/2/3/4$};
\node (13/2/4) at (-2,2)  {$13/2/4$};
\node (1/24/3) at (2,2)  {$1/24/3$};

\draw[thick] (1/2/3/4)--(13/2/4);
\draw[thick] (1/2/3/4)--(1/24/3);

\draw(0,-1.2) node {\large $\NC_G$};

\end{scope}

\end{tikzpicture}
\end{center}
\caption{A graph and its bond lattice and noncrossing bond poset. \label{BL-and-NCBP-Fig}}
\end{figure}
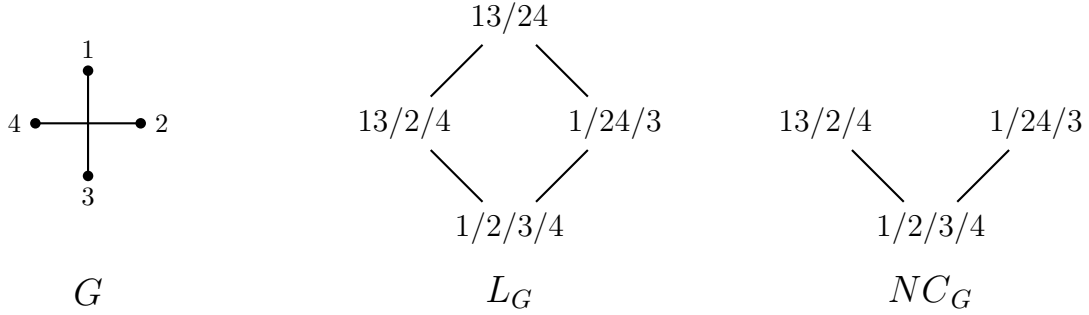

See Figure~\ref{BL-and-NCBP-Fig} for an example of a graph, its bond lattice, and noncrossing bond poset.  Unlike the bond lattice, the noncrossing bond poset is not  always a lattice.  Note that for the graph $G$ in  Figure~\ref{BL-and-NCBP-Fig}, the bond $13/24$ is crossing and so is not in $\NC_G$.  Thus for this graph $G$, $\NC_G$ not only fails to be a lattice, it even fails to have a maximum element.

Note that $\NC_G$ is not necessarily a meet semi-lattice either; the graph $G$ in Figure~\ref{graphFig} gives one such example.  The bonds $H$ and $K$ in Figure~\ref{graphFig} are noncrossing and thus are in $\NC_G$. However, they do not have a meet in $\NC_G$ since the noncrossing bonds $X$ that are below $H$ and $K$ (i.e.~are contained in $H$ and $K$) are those with $E(X) = \emptyset$, $\{14\}$, or $\{35\}$. These do not contain a unique maximal element, so $H$ and $K$ do not have a meet.  

\begin{definition}\label{crossingClosedDef}
Let $G$ be a graph and let $e$ and $f$ be two crossing edges of $G$.  We say $e$ and $f$ are {\em crossing closed} if there exists a unique induced connected subgraph of $G$ containing $e$ and $f$ that is minimal among all such subgraphs with respect to containment.  If such a subgraph exists, we denote it by $J(e,f)$.    We say $G$ is \emph{crossing closed} if every pair of crossing edges in $G$ are crossing closed.
\end{definition}   

Note that the graph $G$ in Figure~\ref{graphFig} is not crossing closed since $14$ and $35$ are crossing but not crossing closed.  There are two distinct minimal, induced, connected subgraphs of $G$ that contain $14$ and $35$, namely the bond $H = G[12345/6]$ and the bond $K = G[13456/2]$ shown in Figure~\ref{graphFig}. On the other hand, trees and complete graphs are crossing closed.  If $G$ is a tree and $e$ and $f$ are two crossing edges, then $J(e,f)$ is the unique path in $G$ with end-edges $e$ and $f$.  If $G$ is a complete graph, then $J(e,f)= G[e \cup f]$.  However, this is far from a complete list of families of crossing closed graphs.  We will save the discussion of more families which are crossing closed  for later sections after we have developed a few more concepts. 

Before we continue, let us explain the choice of the letter ``$J$" in the notation $J(e,f)$.  As we will see in Theorem~\ref{crossingclosededgeThm}, in the case that  $J(e,f)$ exists  it is (essentially) the join of $e$ and $f$ in $\NC_G$. We also wish to explain the reason we use the term ``crossing closed".  As was mentioned in the introduction, the bond lattice is the lattice of flats of the cycle matroid of the graph. One of the many equivalent ways to define a matroid is through a closure operator. In terms of the lattice of flats, the closure of a subset of the ground set is the join of the elements in the lattice of flats.  As we mentioned before, crossing closed edges are exactly the crossing edges which have a join in the graph's noncrossing bond poset and crossing closed graphs are exactly the graphs whose noncrossing bond poset is a lattice.  Thus, a graph being crossing closed implies the existence of a closure operator on the crossing edges that behaves in a similar way that the closure operator does in the cycle matroid.  Note that this does not give us a matroid structure as our closure operator is not the same as the one for matroids.

\begin{lem}\label{mefLemma}
Let $G$ be a graph and let $e$ and $f$ be two crossing closed edges of $G$.  Then $J(e,f)$ is of one of the following two forms, depending on whether or not there is an edge in $G$ connecting a vertex in $e$ to a vertex in $f$.

\begin{enumerate}
    \item There is an edge  in $G$ connecting a vertex of $e$ to a vertex of $f$.  In this case, $J(e,f) = G[e \cup f]$  and {\em $J(e,f)$ is a subgraph of $K_4$}.   
    
    \item There does not exist an edge between a vertex of $e$ and a vertex of $f$.  In this case $J(e,f)$ has the form of the graph depicted in Figure~\ref{dumbellGraphFig}.  Moreover, all vertices in $J(e,f)$ not on $e$ or $f$ are cut vertices of $G$ that separate $e$ and $f$.

\end{enumerate}
\end{lem}

\begin{figure}
\begin{center}
\begin{tikzpicture}[scale=1]
\coordinate (v1) at (0,1);
\coordinate (v2) at (0,-1);
\coordinate (v3) at (6,1);
\coordinate (v4) at (6,-1);

\coordinate (x0) at (1,0);
\coordinate (x1) at (2,0);
\coordinate (xk1) at (4,0);
\coordinate (xk) at (5,0);
\coordinate (dot1) at (2.5,0);
\coordinate (dot2) at (3.5,0);

\draw[thick] (v1)--(v2) (v3) -- (v4);
\draw[dashed] (v1)--(x0)--(v2);
\draw[thick] (x0)--(x1)--(dot1);
\draw[thick] (dot2)--(xk1)--(xk);
\draw[dashed] (v3)--(xk)--(v4);
\node at (3,0) {$\cdots$};

\foreach \v in {v1,v2,v3,v4,x0,x1,xk1,xk} \fill(\v) circle (.1);
\draw(-.4,0) node{\footnotesize $e$};
\draw(-.4,1) node{\footnotesize $v$};
\draw(-.4,-1) node{\footnotesize $v'$};
\draw(6.4,0) node{\footnotesize $f$};
\draw(6.4,1) node{\footnotesize $w$};
\draw(6.4,-1) node{\footnotesize $w'$};
\draw(1,-.5) node{\footnotesize $x_0$};
\draw(2,-.5) node{\footnotesize $x_1$};
\draw(4,-.5) node{\footnotesize $x_{k-1}$};
\draw(5,-.5) node{\footnotesize $x_{k}$};
\draw(3,-2) node{$J(e,f)$};

\end{tikzpicture}
 \end{center}
 \caption {At least one of the dotted edges incident to $x_0$ and at least of the dotted edges incident to $x_k$ exist in $J(e,f)$.}\label{dumbellGraphFig}

\end{figure}
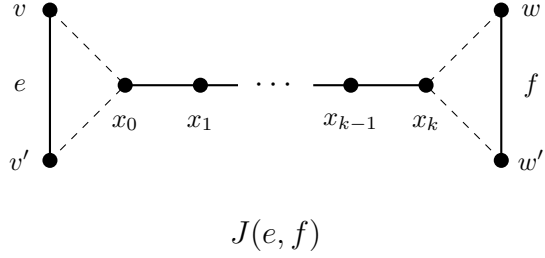
\begin{proof}

Suppose there is an edge  in $G$ connecting a vertex of $e$ to a vertex of $f$.  Then $G[e \cup f]$ is connected, induced, and contains $e$ and $f$.  It is also clearly minimal with respect to those properties so $J(e,f) = G[e \cup f]$.

Now suppose there are no edges  in $G$ connecting a vertex of $e$ to a vertex of $f$.  Let $T$ be a spanning tree of $J(e,f)$ that contains $e$ and $f$. We claim that $T$ cannot have a vertex of degree more than 3 with respect to $T$.  To see why, suppose this was not the case.  Then since $T$ is a tree, it must have at least 3 leaves. Both endpoints of $e$ cannot be leaves of $T$ as then $T$ would then be $e$ and $T$ must also contain $f$. The same can be said of $f$, so $e$ and $f$ together can contain at most two leaves of $T$.  Thus there must be a third leaf $w$ of $T$ with $w$ not in $e$ or $f$.  $J(e,f)\setminus w$ is then induced, connected (since it contains spanning tree $T\setminus w$), and contains $e$ and $f$, a contradiction to the minimality of $J(e,f)$.  Thus $T$ is a path containing $e$ and $f$. The edges $e$ and $f$ must be the two edges at the ends of $T$.  Otherwise at least one end-vertex of $T$, say $w$, would not be on either $e$ or $f$ and $J(e,f)\setminus w$ would again contradict the minimality of $J(e,f)$.

Let the vertices of $J(e,f)$ be given in the order $v,v',x_0, \dots, x_k, w',w$ that they come along the spanning path $T$, where $e = vv'$, $f = ww'$, and $k \geq 0$.  $E(J(e,f))$ is the set of edges these vertices induce. We claim that, besides the edges in $T$, there are no additional edges in $E(J(e,f))$ except possibly for the edges $x_0v$ and $x_kw$ should either of them be present in $G$. 

Suppose, for the sake of obtaining a contradiction, that the claim at the end of the last paragraph is not true. Then $k > 0$.  If $vx_i$ with $i >0$ is present, then $v'vx_ix_{i+1} \dots x_k w'  w$ is a path containing $e$ and $f$ whose vertex set is strictly contained in $V(J(e,f))$ contradicting the minimality of $J(e,f)$.  Similarly, we cannot have an edge $v'x_i$ with $i>0$. By assumption we have no edges from a vertex of $e=vv'$ to a vertex of $f=ww'$. Thus we have no edges from $v$ or $v'$ to any vertex other than $x_0$. Similarly we have no edges from $w$ or $w'$ to any vertex other than $x_k$.  The only remaining possibility is then that we have an edge $x_ix_j$ with $j-i\geq 2$.  But then $vv'x_0\dots x_ix_j\dots x_kw'w$ is a path containing $e$ and $f$ whose vertex set is strictly contained in $J(e,f)$, again contradicting the minimality of $J(e,f)$.

We now claim that each $x_i$ is a cut-vertex of $G$ and separates $e$ and $f$. Suppose instead that $G\setminus x_i$ has a component $C$ that contains $e$ and $f$.  Then by the minimality of $J(e,f)$, $C$ must contain $J(e,f)$ which contains $x_i$, a contradiction.
\end{proof}

\begin{thm} \label{crossingclosededgeThm} Let $G$ be a graph. Let $e$ and $f$ be two crossing edges of $G$.  Then $e$ and $f$ are crossing closed if and only if $e \vee f$ exists in $\NC_G$.  In the  case $e$ and $f$ are crossing closed, $J(e,f)$ is the unique non-trivial component of $e \vee f$ and $e \vee f$ is the bond with edge set $E(J(e,f))$.  Furthermore, $G$ is crossing closed if and only if $e \vee f$ exists for every pair of crossing edges $e$ and $f$.
\end{thm}

\begin{proof}
Suppose $e$ and $f$ are crossing edges which are crossing closed. Then $J(e,f)$ exists. Let $e$ be the bond with edge set $\{e\}$ and $f$ the bond with edge set $\{f\}$.  Let $H$ be the bond with edge set $E(J(e,f))$.  We claim that $H = e \vee f$ in $\NC_G$.  Suppose $H' \in \NC_G$ and we have $e,f \leq H'$ in $\NC_G$.  Edges $e$ and $f$ must be in the same component $C$ of $H'$ otherwise $H'$ would be crossing.  Thus $J(e,f)$ is a subgraph of $C$ and so $H \leq H'$.  Thus $H$ is the unique minimal element of $\NC_G$ that contains bonds $e$ and $f$ and so is $e\vee f$.

Now suppose $H=e \vee f$ exists so that $H \in \NC_G$.  We must have that $e$ and $f$ belong to the same connected component $C$ of $H$, otherwise $H$ would be crossing.  Since $H = e \vee f$ is the minimal element of $\NC_G$ that contains $e$ and $f$, $C$ must be the unique non-trivial component of $H$.  We claim that $J(e,f) = C$.  To show this, let $H'$ be any connected, induced subgraph of $G$ that contains $e$ and $f$. By Definition~\ref{crossingClosedDef}, we must show that $J(e,f)$ is in $H'$.  If we view $H'$ as the bond with edge set $E(H')$ then $H'$ is in $\NC_G$ (since it contains a unique nontrivial connected component and so is noncrossing).  Since $e, f \leq H'$ in $\NC_G$ we must have $H = e \vee f \leq H'$ and thus $C$ is a subgraph of $H'$.

It thus follows that $G$ is crossing closed if and only if $e \vee f$ exists for every pair of crossing edges $e$ and $f$ of $G$.
\end{proof}

\begin{thm}\label{crossingClosedThm}
Let $G$ be a graph. Then $\NC_G$ is a lattice if and only if $G$ is crossing closed. 
Moreover, if $G$ is crossing closed and $H, H' \in \NC_G$, then $H \wedge H' = H \cap H'$.  Thus $\NC_G$ is a meet semi-lattice of $L_G$. 
\end{thm}

\begin{proof}
First, suppose that $\NC_G$ is a lattice.  Then joins exist and in particular the join $e \vee f$ exists for each pair of crossing edges $e,f$ in $G$.  Thus, by the previous theorem, $J(e,f)$ exists for every pair of crossing edges $e,f$ and $G$ is crossing closed.

Now suppose $G$ is crossing closed. Then $G$ must be a noncrossing bond of itself.  If $G$ had two connected components $C_1$ and $C_2$ that cross, then there must be edges $e \in C_1$ and $f \in C_2$ that cross.  There can be no connected subgraph of $G$ that contains both $e$ and $f$, so $J(e,f)$ cannot exist, contradicting the assumption that $G$ was crossing closed.  Thus $G$ is the unique maximum element of $\NC_G$.  Since $\NC_G$ is finite and contains a maximum element, to show $\NC_G$ is a lattice, it suffices to show that it is a meet-semilattice.

We claim that if $H, H'\in\NC_G$, then $H \cap H' \in \NC_G$.  The meet of $H$ and $H'$ in $L_G$ is $H \cap H'$. If  $H \cap H' \notin \NC_G$, there must be crossing edges $e$ and $f$ belonging to different components of $H \cap H'$.  But $J(e,f)$ is a subgraph of both $H$ and $H'$ and so a subgraph of $H \cap H'$, a contradiction.  Thus $H \cap H' \in \NC_G$ and we also have  $H \wedge H' = H \cap H'$ in $\NC_G$.
\end{proof}

In the appendix, we present Algorithm \ref{crossingclosedalg} that decides if $\NC_G$ is a lattice in time polynomial in $n$, the number of vertices of the graph $G$.  It does this by deciding the equivalent question of whether $G$ is crossing closed. We note in the appendix that a brute-force algorithm for this problem could, in the worst case, take time super-exponential in $n$.

A lattice is called \emph{atomic} if every element is the join of a particular subset of atoms, where $\hat{0}$ is considered to be the empty join.  Moreover, we say the lattice is \emph{(upper) semimodular} if whenever $x\wedge y \cover x,y$, then $x,y\cover x\vee y$.   Here  and throughout the paper we use the notation $x \cover y$ to denote that $y$ covers $x$ in the poset. A finite lattice which is both atomic and semimodular  is called \emph{geometric}. 

It is well-known that there is a bijection between geometric lattices and (simple) matroids.  Since the bond lattice of a graph is the lattice of flats for its cycle matroid, every bond lattice is geometric. The situation for the noncrossing bond poset differs.

\begin{prop}\label{latticeProps}
Let $G$ be a crossing closed graph.  Then  we have the following.
\begin{enumerate}
\item[ (a)] $\NC_G$ is atomic.
\item [(b)] $\NC_G$ is a meet-sublattice of $L_G$.
\item [(c)] If $e$ and $f$ cross in $G$, then $J(e,f)$ is the unique non-trivial connected component of $e \vee f$ (the component that contains $e$ and $f$).
\item[(d)] $\NC_G$ is semimodular if and only if $G$ has no crossing edges.   
\item[(e)] $\NC_G$  has a $\hat{1}$.
\end{enumerate}
\end{prop}

Parts (b) and (c) of the previous proposition were proved in Theorem \ref{crossingClosedThm} and Theorem \ref{crossingclosededgeThm} respectively.  The proofs of (a), (d), and (e) are straightforward.

 By Proposition~\ref{latticeProps} part (e), the noncrossing bond poset of any crossing closed graph has a $\hat{1}$. However, this is not the only way to have a maximal element.  For example, the noncrossing bond poset of the graph $G$ in Figure~\ref{graphFig} is not a lattice since $G$ is not crossing closed but still has a $\hat{1}$, which is $G$ itself. The following proposition provides several characterizations for when a $\hat{1}$ exists.  The proof is straightforward and thus omitted.

\begin{prop}\label{1hatProp}
Let $G$ be a graph. Then the following are equivalent. 
\begin{enumerate}
    \item [(a)] $\NC_G$ has a $\hat{1}$.
    \item [(b)] Whenever $e$ and $f$ are crossing edges of $G$, they are in the same connected component of $G$.
    \item [(c)] $G$ is a noncrossing bond of itself.
\end{enumerate}
\end{prop}
 
The following lemma will help us prove results on general graphs by reducing them to the case of connected graphs. 

\begin{lem} \label{lem:product}
Suppose that $G$ consists of connected components $C_1,C_2,\dots, C_k$ such that no edges of $C_i$ and $C_j$ cross for all $i\neq j$. Then $\NC_G \cong NC_{C_1}\times \NC_{C_2}\times \cdots \times \NC_{C_k}$.
\end{lem}

\begin{proof}
Using induction, it suffices to show this result when $k=2$. In that case it is easy to check that the map $\varphi : \NC_G \rightarrow \NC_{C_1}\times \NC_{C_2}$ given by $\varphi(H)  = (H\cap E(C_1), H\cap E(C_2))$ is an isomorphism.
\end{proof}  

\begin{figure}
\begin{center}
\begin{tikzpicture}[scale=.9]
\coordinate (v1) at (0,1);
\coordinate (v2) at (.87,.5);
\coordinate (v3) at (.87,-.5);
\coordinate (v4) at (0,-1);
\coordinate (v5) at (-.87,-.5);
\coordinate (v6) at (-.87,.5);

\draw[thick] (v1)--(v5)--(v3);
\draw[thick] (v1) -- (v4)--(v2)--(v6);

\foreach \v in {v1,v2,v3,v4, v5,v6} \fill(\v) circle (.1);
\draw(0,1.4) node{\footnotesize 1};
\draw (.87+.4,.5+.4) node{\footnotesize 2};
\draw(0,-1.4) node{\footnotesize 4};
\draw (.87+.4,-.5-.4) node{\footnotesize 3};
\draw (-.87-.4,-.5-.4) node{\footnotesize 5};
\draw (-.87-.4,.5+.4) node{\footnotesize 6};
\draw(0,-2.2) node {\large $G$};

\end{tikzpicture}
\end{center}
\caption{A graph whose noncrossing bond poset is not graded. \label{nonGradedFig} }
\end{figure}
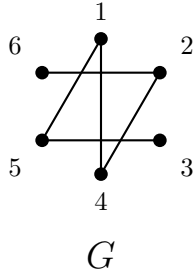

We will now consider the grading of $\NC_G$.  Unlike the bond lattice which is always graded, the noncrossing bond poset need not be graded.  Consider the graph in Figure~\ref{nonGradedFig}.  The bond corresponding to the partition $1/26/35/4$  is noncrossing,  but the only element of $\NC_G$ which covers $1/26/35/4$ is $123456$.   It follows that $1/2/3/4/5/6 \cover 1/26/3/4/5 \cover 1/26/35/4 \cover  123456$ is a maximal chain in $\NC_G$. However, there is another maximal chain  $1/2/3/4/5/6 \cover 14/2/3/5/6 \cover 124/3/5/6 \cover 1246/3/5 \cover 12456/3 \cover 123456$ and so $\NC_G$ is not graded. Note that this graph is a path and hence is crossing closed. Thus $G$ being a crossing closed graph does not imply $\NC_G$ to be graded. 

In the previous example, we were able to find a saturated $\hat{0}$--$\hat{1}$ chain in which each cover is obtained by merging exactly two blocks of the corresponding partition.  It turns out that for graphs whose connected components do not cross (such as connected graphs and crossing closed graphs) such a chain can always be found. In the following proposition and throughout the remainder of the paper, we use $cc(H)$ to denote the number of connected components of a graph $H$.

\begin{prop}\label{gradedProp}
Let $G$ be a graph on $[n]$ which is a noncrossing bond of itself.  Then $\NC_G$ is graded if and only for every cover relation $H\cover H'$, there are exactly two blocks of $H$ that merge to get $H'$.  Moreover, in the case that $\NC_G$ is graded, the rank function is given by $\rho(H)=n-cc(H)$.
\end{prop}

\begin{proof}
 By Proposition~\ref{1hatProp}, $\NC_G$ has a $\hat{1} = G$ and the connected components of $G$ do not cross.  By Lemma \ref{lem:product}, if $G$ has more than one connected component, then $\NC_G$ is the product of noncrossing bond posets, one for each connected component.  Since the product of graded posets is graded and the rank function of the product is the sum of the rank functions, it suffices to prove the result for connected graphs and so we may assume $G$ is connected.

  Let $T$ be a spanning tree of $G$.  Let $e_1,e_2,\dots, e_{n-1}.$ be a sequence of edges of $T$ such that for each forest in the sequence $\{e_1\}, \{e_1,e_2\}, \dots, \{e_1,e_2,\dots, e_{n-1}\}$ there is a unique nontrivial connected component.   For each $1\leq i\leq n-1$, let $H_i$ be the induced subgraph on $\{e_1,e_2,\dots, e_i\}$. Since each $H_i$ has a unique nontrivial connected component, each $H_i$ is noncrossing.  Moreover,  $\hat{0} \cover H_1\cover \cdots \cover H_{n-1} = G$ is a maximal chain of length $n-1$.   Since the minimum element of $\NC_G$ has $n$ blocks, the maximum element has $1$ block, and there is a  maximal chain of length $n-1$, $\NC_G$ is graded if and only if for every cover relation $H\cover H'$, there are exactly two blocks of $H$ that merge to get $H'$.  Now the last statement of the theorem follows immediately since when we merge two blocks, the number of connected components decreases by one.
 \end{proof}

\section{The M\"obius Function and Characteristic Polynomial}\label{MobSec}

In this section, we introduce a family of graphs called upper crossing closed graphs.  The motivation is that this is a class of graphs for which we are able to provide combinatorial interpretations of the M\"obius functions and characteristic polynomials of the corresponding noncrossing bond posets. We briefly recall the definitions of the M\"obius function and characteristic polynomial. For a more information, we refer the reader to~\cite[Chapter 3]{s:ec1}.

We will be dealing with the one-variable version of the M\"obius function,  defined for posets $P$ with $\hat{0}$. This can be recursively defined by
$$
\mu(x) =
\begin{cases}
1 & \mbox{ if } x=\hat{0,}\\
-\displaystyle \sum_{y<x}\mu(y)  & \mbox{ otherwise.}
\end{cases}
$$
Moreover, if $P$ has $\hat{0}$ and also is graded with rank function $\rho$, then  the \emph{characteristic polynomial} of $P$ is given by
$$
\chi(P,t)  = \sum_{x\in P} \mu(x) t^{\rho(P)-\rho(x)}.
$$

Let $\mathrm{ch}(G,t)$ be the chromatic polynomial of the graph $G$, the polynomial $p(t)$ such that for all positive integers $t$, $p(t)$ is the number of proper colorings of $G$ using at most $t$ colors  The chromatic polynomial of a graph $G$ and characteristic polynomial of its bond lattice are related in the following way.
\begin{thm}(see~\cite[Theorem 2.7]{s:iha} ) \label{relationThm}
For all finite graphs $G$, $\mathrm{ch}(G,t) = t^{cc(G)} \chi(L_G,t)$.
\end{thm}

In~\cite{whi:lem}, Whitney gave a combinatorial interpretation for the coefficients of the chromatic polynomial in terms of non-broken circuit sets or NBC sets, which are defined as follows. Let $G$ be a graph.  Let $\unlhd$ be a total ordering on the edges of $G$. A \emph{broken circuit} of $G$ with respect to $\unlhd$ is a collection of edges of $G$ obtained by removing the smallest edge of a cycle of $G$ with respect to that ordering.  We say a subset $S$ of $E(G)$ is a \emph{non-broken circuit set} or {\em NBC set} if $S$ contains no subsets which are broken circuits.  Let $\nbc_k(G)$ be the number of $k$-edge NBC sets of $G$ with respect to $\unlhd$.  Whitney showed the following.

\begin{thm}[Whitney~\cite{whi:lem}]\label{WhiNBCThmchromatic}
Let $G$ be a finite graph on $[n]$.  Then for any total ordering $\unlhd$ of $E(G)$,
$$
\mathrm{ch}(G,t)=\sum_{k\geq 0} (-1)^{k}\nbc_k(G) t^{n-k}.
$$
\end{thm}
Part of the interest of this theorem is its assertion of the non-obvious fact that that the number of NBC sets of size $k$ does not depend on the ordering of the edges.  By Theorem \ref{relationThm} and the fact that $L_G$ has rank function $\rho(H) = n - cc(H)$ we have the following.
\begin{thm}[Whitney~\cite{whi:lem}]\label{WhiNBCThm}
Let $G$ be a finite graph.  Then for any total ordering $\unlhd$ of $E(G)$,
$$
\chi(L_G,t)=\sum_{k\geq 0} (-1)^{k}\nbc_k(G) t^{\rho(L_G)-k}.
$$
\end{thm}

\begin{figure}
\begin{center}
\begin{tikzpicture}
\coordinate (v1) at (0,1);
\coordinate (v2) at (1,0);
\coordinate (v3) at (0,-1);
\coordinate (v4) at (-1,0);
\draw[thick] (v1) -- (v2) -- (v4) --  (v3)--cycle;
\foreach \v in {v1,v2,v3,v4} \fill(\v) circle (.1);
\draw(0,1.4) node{\footnotesize 1};
\draw(0,-1.4) node{\footnotesize 3};
\draw(1.4,0) node{\footnotesize 2};
\draw(-1.4,0) node{\footnotesize 4};
\draw(0,-3.2) node {\large $G$};

\begin{scope}[shift = {(8,-2)}]

\node (1/2/3/4) at (0,0)  {$1/2/3/4$};

\node (12/3/4) at (-5,2)  {$12/3/4$};
\node (13/2/4) at (-2,2)  {$13/2/4$};
\node (1/24/3) at (2,2)  {$1/24/3$};
\node (1/2/34) at (5,2)  {$1/2/34$};

\draw[thick] (1/2/3/4)--(12/3/4);
\draw[thick] (1/2/3/4)--(13/2/4);
\draw[thick] (1/2/3/4)--(1/24/3);
\draw[thick] (1/2/3/4)--(1/2/34);

\node (123/4) at (-5,4)  {$123/4$};
\draw[thick] (12/3/4)--(123/4)--(13/2/4);
\node (124/3) at (-2.5,4)  {$124/3$};
\draw[thick] (12/3/4)--(124/3)--(1/24/3);

\node (12/34) at (0,4)  {$12/34$};
\draw[thick] (12/3/4)--(12/34)--(1/2/34);

\node (134/2) at (2.5,4)  {$134/ 2$};
\draw[thick] (13/2/4)--(134/2)--(1/2/34);

\node (1/234) at (5,4)  {$1/ 234$};
\draw[thick] (1/24/3)--(1/234)--(1/2/34);

\node (1234) at (0,6)  {$1234$};
\draw[thick] (123/4)--(1234)--(124/3);
\draw[thick] (12/34)--(1234)--(134/2);
\draw[thick] (1/234)--(1234);

\draw(0,-1.2) node {\large $\NC_G$};

\tikzstyle{every node}=[red, inner sep=0.4pt, scale=1, minimum width=4pt]

\draw (1.25, 0) node {$+1$};

\draw (-3.75, 2) node {$-1$};
\draw (-.75, 2) node {$-1$};
\draw (3.25, 2) node {$-1$};
\draw (6, 2) node {$-1$};

\draw (-4, 4) node {$+1$};
\draw (-1.5, 4) node {$+1$};
\draw (1, 4) node {$+1$};
\draw (3.5, 4) node {$+1$};
\draw (6, 4) node {$+1$};

\draw (1.25, 6) node {$-2$};

\end{scope}
\end{tikzpicture}
\end{center}
\caption{The twisted 4-cycle and its noncrossing bond poset. M\"obius values are in red. \label{twistedC4Fig} }
\end{figure}
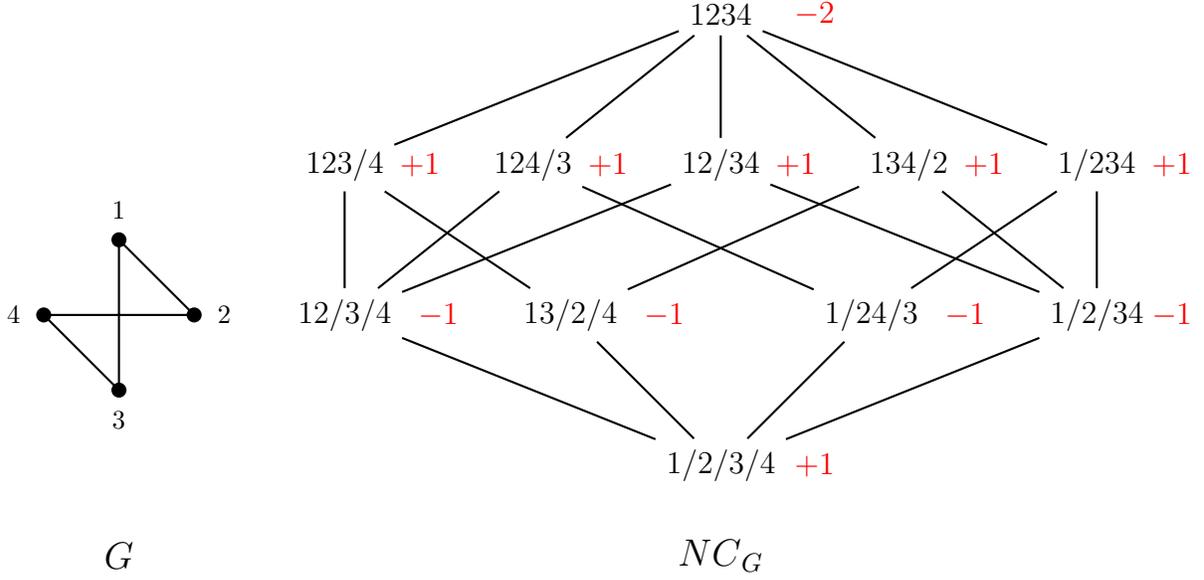

As an example,  consider the twisted 4-cycle graph $G$ in Figure~\ref{twistedC4Fig}.  It is not hard to calculate that
$$\chi(L_G,t) = t^3 -4t^2+ 6t -3.$$
This can be done by calculating the M\"obius function of $L_G$ and then using the definition of the characteristic polynomial or by calculating $\mathrm{ch}(G,t)$ and then using Theorem \ref{relationThm}. We can also use Theorem \ref{WhiNBCThm}.  Suppose $\unlhd$ is lexicographic order, i.e.~$12 \lhd 13 \lhd 24 \lhd 34$. Since $G$ is a cycle, the only broken circuit is $\{13, 24,34\}$ and hence every subset of $E(G)$ is an NBC set except for $\{13, 24, 34\}$ and $\{12,13, 24,34\}$.  Thus the absolute values of the coefficients of $\chi(L_G,t)$ are indeed $\nbc_0(G) =1, \nbc_{1}(G)=4, \nbc_{2}(G)=6, \nbc_{3}(G)=3$ and $\nbc_k(G)=0$ for $k>3$.  This agrees with the coefficients one finds by using the M\"obius function or the chromatic polynomial.

Now let us compare this with the characteristic polynomial of $\NC_G$.  From the M\"obius values shown in  Figure~\ref{twistedC4Fig}, we see that 
$$
\chi(\NC_G,t ) = t^3- 4t^2 +5t -2.
$$
Since the absolute value of the coefficients of  $\chi(\NC_G,t)$ are less than or equal to the corresponding values in $\chi(L_G,t)$, it is plausible that the coefficients of $\chi(\NC_G,t)$ might count a subset of the NBC sets of $G$. 

If $\unlhd$ is a total order of $E(G)$ we say that $S$ is a \emph{noncrossing NBC set} with respect to $\unlhd$, if $S$ is an NBC set with respect to that order and contains no edges which cross in the graphical representation of $S$. We define $\NCNBC_k(G)$ to be the set of noncrossing NBC sets of size $k$ and $\ncnbc_k(G)$ to be the number of such sets.  

If $\unlhd$ is again the lexicographic order on the edges of the twisted 4-cycle, then all NBC sets are also noncrossing NBC sets except for  $\{13, 24\}$ and $\{12, 13, 24\}$.  Thus the sequence $\ncnbc_0(G)=1,\ncnbc_1(G)=4, \ncnbc_2(G) = 5, \ncnbc_3(G)=2$ and $\ncnbc_k(G)=0$ for $k>3$ does match the sequence of absolute values of the coefficients of $\chi(\NC_G,t)$.  Unfortunately, unlike $\nbc_k(G)$, $\ncnbc_k(G)$ does depend on the ordering $\unlhd$ (despite the fact that the notation does not reflect this).  For example, if we order the edges as $13 \lhd 24 \lhd 12 \lhd 34$, there is only one noncrossing NBC set with $3$ edges, instead of $2$ and so this ordering will not give the correct coefficients. However, the upper crossing closed graphs (defined below), do have edge orderings for which the M\"obius values of $\NC_G$ and the absolute values of the coefficients of $\chi(\NC_G,t)$ are indeed the counts of noncrossing NBC sets with respect to those orderings.

\begin{definition}\label{upperCrossingClosedDef}
We say a graph $G$ is {\em upper crossing closed} if it is crossing closed and there is a total ordering $\unlhd$ on $E(G)$ such that for every pair of crossing edges $e$ and $f$, $J(e,f)$ contains an edge $h$ such that $h \lhd e,f$.  If $\unlhd$ is one such ordering, we say that $G$ is {\em upper crossing closed with respect to $\unlhd$} and also that {\em $\unlhd$ is an upper crossing closed ordering of $E(G)$}.
\end{definition}

Note that the twisted 4-cycle in Figure~\ref{twistedC4Fig} is upper crossing closed with respect to lexicographic order since $J(13,24) = G$ contains $12$ and $12 \lhd 13, 24$.  Note also that it is not upper crossing closed with respect to the other ordering $13\lhd 24 \lhd 12 \lhd 34$ we considered previously since $13$ is then the minimum edge in $J(13,24)$.  As we saw earlier, when using the lexicographic ordering on the twisted 4-cycle (which is upper crossing closed),  the coefficients of the characteristic polynomial count noncrossing NBC sets.   This is not a coincidence as we see next.

\begin{thm}\label{NCNBCThm}
Let $G$ be a graph on $[n]$. If $G$ is  upper crossing closed  with respect to the order $\unlhd$ on $E(G)$.  Then for all $H \in \NC_G$, 
$$
\mu(H) = (-1)^{n-cc(H) } \ncnbc_{n-cc(H)}(H).
$$
Moreover, if $\NC_G$ is graded, then
$$
\chi(\NC_G, t) =\sum_{k\geq 0} (-1)^k\ncnbc_k(G) t^{\rho(\NC_G) -k}.
$$

\end{thm}

Before we can prove Theorem \ref{NCNBCThm}, we need to discuss the notion of NBB sets developed by Blass and Sagan in~\cite{bs:mfl}.

\begin{definition} [Bass and Sagan \cite{bs:mfl}]
Let $L$ be a lattice and let $\unlhd$ be a partial order on the atoms of $L$.   A subset $S$ of the atoms of $L$ is \emph{bounded below} if there exists an atom $a$ such that
\begin{enumerate}
\item[(a)]  $a \lhd s$  for all $s\in S$
\item[(b)] $a < \bigvee S$
\end{enumerate}
We say a subset $S$ of the atoms of $L$ is a \emph{non-bounded below (NBB)} set for $x$ if $S$ contains no bounded below sets and $\bigvee S = x$.  
\end{definition}

Blass and Sagan's result generalizes Whitney's NBC theorem. In particular, if one considers a graph $G$ and its bond lattice $L_G$, then the edges of $G$ correspond to the atoms of $L_G$ and the NBB sets are exactly the NBC sets of the graph with respect to whatever ordering is put on the atoms/edges. Blass and Sagan showed that we can use NBB sets to determine the value of the M\"obius function.

\begin{thm}[Blass and Sagan \cite{bs:mfl}]
Let $L$ be a lattice and let  $\unlhd$ be a partial order on the atoms of $L$. Then for all $x\in L$, 
$$
\mu(x) = \sum_{B} (-1)^{|B|}  
$$
where the sum is over NBB sets for $x$.
\end{thm}

In Lemma \ref{NBBisNCBNC} below, we will show that if $\unlhd$ is an upper crossing closed ordering of a graph $G$ then the NBB sets and the NCNBC sets are the same.  We illustrate this using our running example of the twisted 4-cycle in Figure~\ref{twistedC4Fig}. For this example, we use the lexicographic ordering which was already shown to be upper crossing closed.  First, let us note that the empty set and any singleton subset of atoms is  NBB since their joins only have themselves below them.  Moreover, since $G$ is a 4-cycle, any subset of 2 edges which do not cross is NBB for the same reason.  Now let us turn our attention to the two edges that do cross, $13$ and $24$.   Their join is the entire graph and since $12$ is lexicographically smaller than $13$ and $24$,  $\{13,24\}$ is a bounded below set.   So every 2-element subset of the atoms except $\{13,24\}$ is  NBB.  Finally, we consider the subsets of size 3.  Since $G$ is a 4-cycle the join of any 3-element subset is the entire graph. Since $12$ is the smallest edge lexicographically, $\{13,24,34\}$ is bounded below.  Of the remaining 3-element subsets,  $\{12,13,24\}$ is not NBB since it contains $\{13,24\}$, $\{12,13,34\}$ is NBB, and $\{12,24,34\}$ is NBB.  We conclude that the NBB sets of the twisted 4-cycle with edges ordered lexicographically are $\emptyset$, $\{12\}$, $\{13\}$, $\{24\}$, $\{34\}$, $\{12, 13\}$, $\{12, 24\}$, $\{12,  34\}$, $\{13,34\}$, $\{24,34\}$, $\{12,13, 34\}$, $\{12, 24, 34\}$, which match up exactly with the noncrossing NBC sets of the twisted 4-cycle for the lexicographic ordering.

\begin{lem}\label{sameJoin}
Let $G$ be a  crossing closed  graph and let $S\in\NCNBC_k(G)$.  Then the join of the elements in $S$ is the same in $L_G$ and $\NC_G$
\end{lem}
\begin{proof}
In this proof, we will use $\bigvee_{L_G}$ and $\bigvee_{\NC_G}$ to denote the join operators in $L_G$ and $\NC_G$ respectively. First, let us show that $\bigvee_{L_G} S$ is a noncrossing bond.  Suppose this was not the case and let $C_1$ and $C_2$ be components of $\bigvee_{L_G} S$ which have crossing edges. Let $S_1=S\cap E(C_1)$ and $S_2= S\cap E(C_2)$.  Then $S_1$ is spanning tree of $C_1$ and $S_2$ is a spanning tree of $C_2$.  Since $C_1$ and $C_2$ cross, there exists edges $ac\in E(C_1)$ and $bd\in E(C_2)$ with $a<b<c<d$ which cross.  In $S_1$ there is a path from $a$ to $c$, but this path must separate $b$ and $d$.  Similarly, in $S_2$ there is a path between $b$ and $d$.  This path must cross the path between $a$ and $c$. However, these two paths are in different connected components. This implies that $S_1$ and $S_2$ must cross,  but then $S$ is crossing which is impossible.

Since $S$ is a collection of edges of $G$ and contains no broken circuits, it forms a spanning forest of $G$.  It is not hard to see that $\bigvee_{L_G} S$ is the bond whose induced connected components are the connected components of $S$.  As we saw, since $S$ is noncrossing,  $\bigvee_{L_G} S$ is a noncrossing bond in $L_G$ and hence is an element of $\NC_G$.   It follows that  the partition associated with $\bigvee_{L_G} S$ is noncrossing.  It is not hard to see that this is exactly the same partition associated with  $\bigvee_{\NC_G} S$.  Thus, the result holds.
\end{proof}

\begin{lem}\label{NBBisNCBNC}
Let $G$ be an upper crossing closed graph with total ordering $\unlhd$ on $E(G)$.  Suppose $G$ is  upper crossing closed  with respect to $\unlhd$. Order the atoms of $\NC_G$ by $\unlhd$. Then $\NBB_k(G)= \NCNBC_k(G)$, where $\NBB_k(G)$ is the set of non-bounded below sets of $\NC_G$ with $k$ elements.
\end{lem}
\begin{proof}
Suppose that $S\in\NBB_k(G)$, but that $S\notin \NCNBC_k(G)$.  If $S$ is not an NBC set, then it contains a broken circuit $C$.  Let $e$ be the edge removed from the cycle to obtain $C$.  Then it is not hard to see that $C$ is a bounded below set with $e$ as the atom which is below all the elements of $C$.  This would imply that $S\notin\NBB_k(G)$.  Thus, $S$ must be an NBC set.  Now suppose that $S$ has crossing edges.  Let $S'$ be a set consisting of two such crossing edges.  Since $G$ is upper crossing closed, $\bigvee S'$ contains an edge smaller than all the edges of $S'$.  It follows that $S'$ is a bounded below set, but then $S$ is not an NBB set as it contains $S'$. Thus $S\in \NCNBC_k(G)$.

Next, suppose that $S\in \NCNBC_k(G)$, but that $S\notin \NBB_k(G)$.  Then $S$ contains a bounded below set, $T$. Note that  since $T\subseteq S$, $T$ is a noncrossing NBC set. Moreover, since $T$ is bounded below, there exists an atom $e$  such that $e\lhd t$ for all $t\in T$ and $e<\bigvee T$.  By Lemma~\ref{sameJoin}, $\bigvee T$ is the same in $L_G$ and $\NC_G$. Thus, $e< \bigvee T$ and $e\notin T$, implies that $e$ must be in some  cycle $C$ of $\bigvee T$.  To see why, note that $\bigvee T$ is the induced subgraph on $T$.  The fact that $e< \bigvee T$, implies $e$ is in $\bigvee T$.  So if $e=uv$, there is a path from $u$ to $v$ in $\bigvee T$.  The set $T$ must contain a spanning tree for the component containing $u$ and $v$.  But since $e\notin T$, $e$ is not on this spanning tree, so $e$ must be on a cycle.  Moreover, since $e$ is smaller than all the elements of $T$, $C \setminus e$ would be a broken circuit of $T$. But then $S$ is not an NBC set which is impossible.
\end{proof}

We are now ready to prove Theorem~\ref{NCNBCThm}.

\begin{proof}[(Proof  of Theorem~\ref{NCNBCThm})]    Order the atoms of $\NC_G$ by $\unlhd$.   Using the fact that $G$ is upper crossing closed, Lemma~\ref{NBBisNCBNC} shows that a subset of atoms of $\NC_G$ is NBB if and only if it is  a noncrossing NBC set of $G$.  Then using Blass and Sagan's result, we have that for each $H \in \NC_G$,
$$
\mu(H) =  \sum_{B} (-1)^{|B|} 
$$
where the sum is over all the noncrossing NBC sets $B$ such that $\bigvee B = H$.  Since $B$ is a noncrossing NBC set, Lemma~\ref{sameJoin} implies that  $\bigvee B$ is the same in $L_G$ and $\NC_G$. We claim that for a fixed $H$, all the NBC sets whose join is $H$ in $L_G$ have the same size, namely $n-cc(H)$. To see why, suppose that  $S$ is an NBC set and $\bigvee S =H$.  It must be the case that the edges in $S$ form a subgraph of $G$ so that its connected components are exactly the connected components of $H$.  Moreover, since NBC sets cannot contain cycles, $S$ must be minimal with respect to spanning the connected components.  So each connected component of $S$ must be a tree.   Thus the number of edges in $S$ is $n-cc(S) = n-cc(H)$. It now follows that 
\begin{align*}
\mu(H)& =  \sum_{B} (-1)^{|B|}\\
  &= \sum_{B} (-1)^{n-cc(H) }\\
&=  (-1)^{n-cc(H) }  \# (\mbox{noncrossing NBC sets of $G$ whose join is $H$})\\
&=(-1)^{n-cc(H) } \ncnbc_{n-cc(H)}(H).
\end{align*}

To finish, note that if $\NC_G$ is graded, then since $G$ is crossing closed, it has a $\hat{1}$ and so Proposition~\ref{gradedProp} implies that the rank function  of $\NC_G$ is $\rho(H)= n-cc(H)$.  Since $\chi(\NC_G,t) = \sum_{H\in \NC_G} \mu(H)t^{\rho(\NC_G)-\rho(H)}$,
\begin{align*}
\chi(\NC_G,t) &=\sum_{H\in \NC_G} (-1)^{n-cc(H) } \# (\mbox{noncrossing NBC sets of $G$ with join $H$})t^{\rho(\NC_G)-\rho(H)}\\
&= \sum_{k\geq 0} \sum_{\rho(H)=k}  (-1)^{k} \#(\mbox{noncrossing NBC sets of $G$ with join $H$})t^{\rho(\NC_G)-k}\\
&=\sum_{k\geq 0} (-1)^{k} \ncnbc_k(G)t^{\rho(\NC_G)-k}
\end{align*}
as claimed.
\end{proof}

\begin{figure}[H]
\begin{center}
\begin{tikzpicture}[scale=.9]
\coordinate (v1) at (0,1);
\coordinate (v2) at (.87,.5);
\coordinate (v3) at (.87,-.5);
\coordinate (v4) at (0,-1);
\coordinate (v5) at (-.87,-.5);
\coordinate (v6) at (-.87,.5);

\draw[thick] (v1)--(v5)--(v3);
\draw[thick] (v1) -- (v4)--(v2)--(v6);

\foreach \v in {v1,v2,v3,v4, v5,v6} \fill(\v) circle (.1);
\draw(0,1.4) node{\footnotesize 1};
\draw (.87+.4,.5+.4) node{\footnotesize 2};
\draw(0,-1.4) node{\footnotesize 4};
\draw (.87+.4,-.5-.4) node{\footnotesize 3};
\draw (-.87-.4,-.5-.4) node{\footnotesize 5};
\draw (-.87-.4,.5+.4) node{\footnotesize 6};
\draw(0,-2.2) node {\large $G$};

\begin{scope}[shift={(5,0)}]

\coordinate (v1) at (0,1);
\coordinate (v2) at (.87,.5);
\coordinate (v2') at (1,0);
\coordinate (v3) at (.87,-.5);
\coordinate (v4) at (0,-1);
\coordinate (v5) at (-.87,-.5);
\coordinate (v6) at (-.87,.5);
\coordinate (v6') at (-.5,.87);

\draw[thick] (v1)--(v6')--(v5)--(v3);
\draw[thick] (v1) -- (v4)--(v2')--(v2)--(v6);

\foreach \v in {v1,v2,v2',v3,v4, v5,v6,v6'} \fill(\v) circle (.1);
\draw(0,1.4) node{\footnotesize 1};
\draw (.87+.4,.5+.4) node{\footnotesize 2};
\draw (1+.4,0) node{\footnotesize 2$'$};
\draw(0,-1.4) node{\footnotesize 4};
\draw (.87+.4,-.5-.4) node{\footnotesize 3};
\draw (-.87-.4,-.5-.4) node{\footnotesize 5};
\draw (-.87-.4,.5+.4) node{\footnotesize 6};
\draw (-.5-.4,.87+.4) node{\footnotesize 6$'$};
\draw(0,-2.2) node {\large $H$};

\end{scope}

\end{tikzpicture}
\end{center}
\caption{$\NC_G$ and $\NC_H$ are not ranked.  $H$ is upper crossing closed ($G$ is not).
\label{nonGradedUpperCCFig} }
\end{figure}
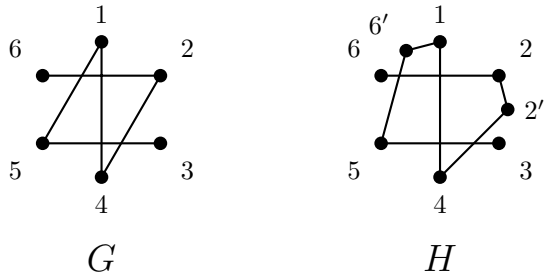

In the statement of Theorem~\ref{NCNBCThm}, we had to make the assumption that $\NC_G$ is graded in order to describe the characteristic polynomial.  We must do so because there exists upper crossing closed graphs with the property that their noncrossing bond posets are not graded.  To construct such an example, we will take a subdivision of the graph $G$ given in Figure~\ref{nonGradedUpperCCFig}. As mentioned earlier, the noncrossing bond poset of the graph  of $G$ is not graded. Additionally, since every edge  of the graph is crossed, it cannot be upper crossing closed.  This is because if $e$ is the smallest edge and it crosses some edge $f$, $J(e,f)$ would need to have an edge smaller than $e$.

Consider the graph $H$ given in Figure~\ref{nonGradedUpperCCFig} which is obtained from subdividing the edges $24$ and $15$ and labeling the new vertices $2'$ and $6'$.  $H$ is also a tree and hence is crossing closed.  Any ordering of the edges in which $16'$ and $22'$ come first is an upper crossing closed ordering since every $J(e,f)$ contains one of these edges.   To see why $\NC_H$ is not graded, consider the bond corresponding to the partition $16'/22'6/35/4$.  If we try to add any edge to this bond, we create a new crossing.  Since $H$ is a tree, if $\NC_H$ was graded, each covering relation would be obtained by adding a single edge, a contradiction.

\begin{figure}[H] 
\begin{center}
\begin{tikzpicture}[scale=.9]

\coordinate (v1) at (0,1);
\coordinate (v2) at (.88,.47);
\coordinate (v3) at (.72,-.69);
\coordinate (v4) at (-.72,-.69);
\coordinate (v5) at (-.88,.47);
\draw[thick] (v1)--(v3)--(v5)--(v2)--(v4)--(v1);
\foreach \v in {v1,v2,v3,v4,v5} \fill(\v) circle (.1);
\draw(0,1.4) node{\footnotesize 1};
\draw(.88+.4,.47) node{\footnotesize 2};
\draw(.72+.4,-.69) node{\footnotesize 3};
\draw(-.72-.4,-.69) node{\footnotesize 4};
\draw(-.88-.4,.47) node{\footnotesize 5};

\draw(0,-2) node {$G$};

\end{tikzpicture}
\end{center}
\caption{5-pointed star \label{5pointStarFig} }
\end{figure}
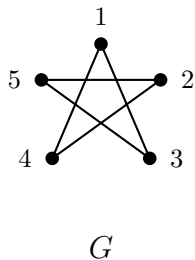

In this section, we saw that the notion of upper crossing closed allows us to give a combinatorial description of the M\"obius function and characteristic polynomial in terms of noncrossing NBC sets.  As we will see in a later section, the assumption of upper crossing closed is not always necessary.  Indeed there are graphs, which are not even crossing closed which have such an interpretation.  However, there are also graphs for which no such interpretation is possible.  Consider the 5-pointed star in Figure~\ref{5pointStarFig}. The characteristic polynomial is given by 
$$
\chi(\NC_G,t) = t^4 - 5t^3 + 5t^2 - 1.
$$
The coefficient of $t$ is 0, but the coefficient of $t^0$ is nonzero.  If the characteristic polynomial was the generating function for noncrossing NBC sets, then this would imply that there is 1 noncrossing NBC set of size $4$, but none of size $3$.  However, every subset of a noncrossing NBC set is a noncrossing NBC set so this is impossible.  This argument generalizes.  Say a polynomial $P(t) = c_nt^n+c_{n-1}t^{n-1}+\cdots +c_1t+c_0$ has an \emph{internal zero} if there exists a $k$ such that $c_k=0$, but $c_{k+1},c_{k-1} \neq 0$.  By the hereditary property of noncrossing NBC sets, the generating function for noncrossing NBC sets cannot have internal zeros. 

We give an algorithm, Algorithm \ref{uppercrossingclosedalg} in the appendix, that when given a graph $G$ will either produce an specific upper crossing closed ordering of $E(G)$ or will produce what we term an obstruction, a specific subgraph of $G$ that clearly shows there can be no such ordering. This gives a forbidden subgraph characterization of upper crossing closed graphs, Theorem \ref{ucccharacterization}: a graph is upper crossing closed if and only if it contains no such subgraph. 

We also prove that Algorithm \ref{uppercrossingclosedalg} will run in time polynomial in $n$, the number of vertices of $G$ and note that a brute-force algorithm will, in the worst case, take time super-exponential in $n$.

\section{Edge Labelings and Shellability}\label{edgeLabSec}

In this section we discuss edge labelings of the noncrossing bond poset. We pay particular attention to the minimum labeling, introduced by Bj\"orner, which is an EL-labeling for every geometric lattice. While the minimum labeling is not  an EL-labeling for every noncrossing bond poset, we give a sufficient condition which guarantees it is an EL-labeling.  We do this for two reasons.  First, this will show that the poset is shellable.  Second, it will allow us to show that the combinatorial interpretation for the M\"obius function in terms of noncrossing NBC sets holds for more than just upper crossing closed graphs.

We briefly review edge labelings of posets.  We refer the reader to~\cite{w:pt} for more information.  Let $P$ be a graded poset.  An \emph{edge labeling} of $P$ is a function $\lambda: \mathcal{E}(P) \rightarrow \Lambda $ where $\mathcal{E}(P)$ is the set of edges of the Hasse diagram of $P$ and $\Lambda$ is a set of labels which is partially ordered.    We note here that although the labels are allowed to be partially ordered, in this article they will always be totally ordered.  Now suppose that $P$ is a graded poset with edge labeling $\lambda$.   Let $\c: x_0\cover x_1 \cover \cdots \cover x_k$ be a saturated chain in $P$.  We say $\c$ is \emph{increasing} if 
$$
\lambda(x_0\cover x_1)<\lambda(x_1\cover x_2)<\cdots< \lambda(x_{k-1}\cover  x_k).
$$
  Moreover, we say $\c$ is \emph{decreasing} if 
$$
\lambda(x_0 \cover  x_1) \geq \lambda(x_1 \cover x_2) \geq \cdots\geq  \lambda(x_{k-1}\cover x_k).
$$
\noindent Before we move on, let us note that while the inequalities for decreasing are allowed to be weak, in this paper, they are always strict.  Thus, there is no need to add the adjective ``weak" to decreasing.

Let $\lambda$ be an edge labeling of $P$. We say $\lambda$ is an \emph{EL-labeling} if every interval has a unique  increasing maximal chain and this chain precedes every other maximal chain in the interval in lexicographic order. Bj\"orner~\cite{b:scmpos} and  Bj\"orner and Wachs~\cite{bw:lsp} showed that there are several nice topological consequences of a poset having an EL-labeling.  For example, they showed that given a poset with an EL-labeling, the order complex of $P$ is shellable and has the homotopy type of a wedge of spheres.  Because of this connection with the topology of the order complex, EL-labelings also have implications for the M\"obius function. In particular, we have the following simple combinatorial interpretation for the M\"obius function for graded EL-labeled posets. 

\begin{thm}[Bj\"orner~\cite{b:scmpos}]\label{mobFuncEL}
Let $P$ be a graded poset with an EL-labeling.  Then
$$
\mu(x) = (-1)^{\rho(x)}\#(\mbox{decreasing saturated chains from $\hat{0}$ to $x$})
$$
\end{thm}

We will now consider an edge labeling for the bond lattice that we can also can apply to the noncrossing bond poset.

\begin{definition}[Bj\"orner~\cite{b:scmpos}]\label{minELLabeling}
Let $G$ be a graph.  Fix some total order $\unlhd$ on $E(G)$. The \emph{minimum labeling} of $L_G$ is defined by
$$
\lambda(H\cover H')  = \min ( E(H')\setminus E(H))
$$
where the minimum is taken with respect to $\unlhd$.
\end{definition}

\begin{figure}
\begin{center}
\begin{tikzpicture}[scale=.7]
\coordinate (v1) at (0,1);
\coordinate (v2) at (1,0);
\coordinate (v3) at (0,-1);
\coordinate (v4) at (-1,0);
\draw[thick] (v3)--(v1) -- (v2) -- (v4);
\foreach \v in {v1,v2,v3,v4} \fill(\v) circle (.1);
\draw(0,1.4) node{\footnotesize 1};
\draw(0,-1.4) node{\footnotesize 3};
\draw(1.4,0) node{\footnotesize 2};
\draw(-1.4,0) node{\footnotesize 4};
\draw(0,-3.2) node {\large $G$};

\begin{scope}[shift = {(8,-2)}]

\node (1/2/3/4) at (0,0)  {$1/2/3/4$};

\node (12/3/4) at (0,2)  {$12/3/4$};
\node (13/2/4) at (-5,2)  {$13/2/4$};
\node (1/24/3) at (5,2)  {$1/24/3$};

\draw[thick] (1/2/3/4)--(12/3/4);
\draw[thick] (1/2/3/4)--(13/2/4);
\draw[thick] (1/2/3/4)--(1/24/3);

\node (123/4) at (-5,4)  {$123/4$};
\draw[thick] (12/3/4)--(123/4)--(13/2/4);
\node (124/3) at (5,4)  {$124/3$};
\draw[thick] (12/3/4)--(124/3)--(1/24/3);

\node (1234) at (0,6)  {$1234$};
\draw[thick] (123/4)--(1234)--(124/3);
\draw(0,-1.2) node {\large $\NC_G$};
\tikzstyle{every node}=[blue, scale=1.1]
\node at (-2.5, 1.45) {$13$};
\node at (.6, 1.2) {$12$};
\node at (2.5, 1.45) {$24$};
\node at (-2.5, 3.45) {$13$};
\node at (2.5, 3.45) {$24$};
\node at (5.5, 3.0) {$12$};
\node at (-5.5, 3.0) {$12$};
\node at (-2.5, 5.45) {$24$};
\node at (2.5, 5.45) {$13$};
\end{scope}
\end{tikzpicture}
\end{center}
\caption{A graph and its noncrossing bond posets labeled by the minimum labeling\label{minLabelFig} }
\end{figure}
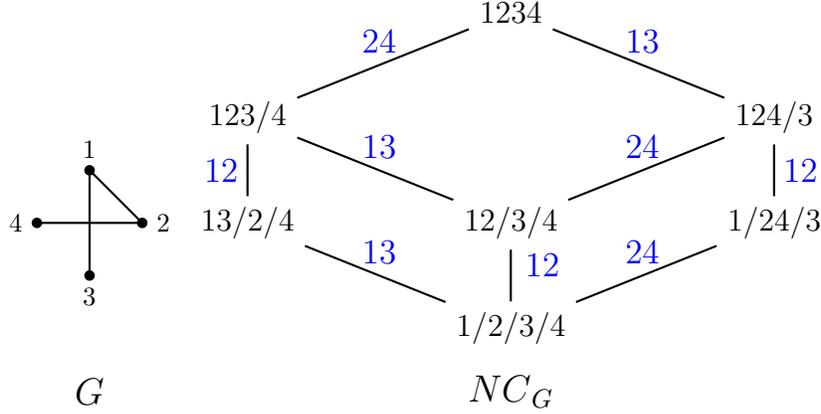

Figure~\ref{minLabelFig} contains an example of the minimum labeling where the edges are ordered lexicographically.  We note that the usual definition of the minimum labeling is phrased in terms of the join of elements. This  makes sense for the bond lattice, however, because the noncrossing bond poset need not be a lattice, we have phrased it in a different (but equivalent) way.  Bj\"orner showed the following concerning the minimum labeling of the bond lattice.  In the following theorem and throughout the section, we will use the term \emph{spanning NBC set} of  a graph $X$ to mean an NBC set $S$ of $X$ such that the induced subgraph $S$ is $X$.

\begin{thm}[Bj\"orner~\cite{b:scmpos}]\label{minLabelNBC}
Let $G$ be a graph and let $\unlhd$ be a total ordering of $E(G)$. Then we have the following (where the NBC sets are taken with respect to $\unlhd$).
\begin{enumerate}
     \item [(a)] The minimum labeling with respect to $\unlhd$ is an EL-labeling of $L_G$ and so $L_G$ is shellable.
    \item [(b)] The labels along any decreasing saturated chain from $\hat{0}$ to $X$ form a spanning NBC set of $X$. 
    \item [(c)] For every $X\in L_G$, each spanning NBC set of $X$ appears exactly once as a saturated decreasing chain from $\hat{0}$ to $X$.
    \end{enumerate}
\end{thm}

Note that Theorem~\ref{mobFuncEL} and Theorem~\ref{minLabelNBC} together imply Whitney's theorem (Theorem~\ref{WhiNBCThm}) concerning the NBC set interpretation of the M\"obius function.   Even though the minimum labeling is an EL-labeling for the bond lattice of any graph, it need not be an EL-labeling for the noncrossing bond poset.  In the next theorem, we give sufficient conditions for the labeling to be an EL-labeling.  We also show that in this setting, we get noncrossing analogues of Bj\"orner's result on NBC sets described in the previous theorem.  As a result,  we get the same combinatorial interpretation for the M\"obius function and characteristic polynomial as we did with upper crossing closed graphs, see Theorem \ref{NCNBCThm}.

\begin{thm}\label{minLabelNC}
Let $G$ be a graph on $[n]$ such that $\NC_G$ has a $\hat{1}$. Let $\unlhd$ be a total ordering of $E(G)$.  Suppose that whenever $H< H'$ and $e=\min E(H') \setminus E(H)$, the bond induced on $E(H) \cup \{e\}$ is noncrossing.  Then we have the following (where the noncrossing NBC sets are taken with respect to $\unlhd$).
\begin{enumerate}
    \item[(a)] $\NC_G$ is graded and $\rho(X) = n -cc(X)$.
    \item [(b)] The minimum labeling with respect to $\unlhd$ is an EL-labeling and so $\NC_G$ is shellable.
    \item [(c)] The labels along any decreasing saturated chain from $\hat{0}$ to $X$ form a spanning noncrossing NBC set of $X$. 
    \item [(d)] For every $X\in \NC_G$, each spanning noncrossing NBC set of $X$ appears exactly once as a saturated decreasing chain from $\hat{0}$ to $X$.
    \item[(e)]  For $H \in \NC_G$, 
$$
\mu(H) = (-1)^{n-cc(H) } \ncnbc_{n-cc(H)}(H)
$$
and
$$
\chi(\NC_G, t) =\sum_{k\geq 0} (-1)^k\ncnbc_k(G) t^{\rho(\NC_G) -k}.
$$
\end{enumerate}
\end{thm}
\begin{proof}
For (a), suppose that $H\cover H'$.  Let $e=\min E(H')\setminus E(H)$ and let $H''$ be the bond induced by $E(H)\cup \{e\}$.  By assumption $H''$ is noncrossing.  Thus, in $\NC_G$ we have $H<H''\leq H'$.  It follows that $H''=H'$ and so by Proposition~\ref{gradedProp}, $\NC_G$ is graded and $\rho(X) = n-cc(X)$.

Now let us show (b). Let $\lambda$ be the minimum labeling with respect to $\unlhd$.  Since $\NC_G$ is a graded subposet of $L_G$ with the same rank function, the sequence of labels along any saturated chain in $\NC_G$ also appears as the labels along that chain in $L_G$.  Since $\lambda$ is an EL-labeling of $L_G$, every interval has at most one increasing maximal chain and  this chain lexicographically precedes all other chains in that interval.   Thus it suffices to show that each interval of $\NC_G$ has an increasing maximal chain.  

Consider an interval $[H,H']$ in $\NC_G$.  We will induct on the length of  $[H,H']$.  If the length is 1, the result is trivial. Now assume that the length is larger than $1$.  Let $e=\min E(H')\setminus E(H)$ and let $H''$ be the bond induced on $E(H)\cup \{e\}$.  By assumption $H''$ is noncrossing and so is in $[H,H']$.  Then $\lambda(H\cover H'')=e$.  By the inductive hypothesis $[H'',H']$ has an increasing maximal chain which starts with a label larger than $e$.  Concatenating this chain with $H\cover H''$ will produce an increasing maximal chain of $[H,H']$.  It follows that the minimum labeling is an EL-labeling and so $\NC_G$ is shellable.

Next we show part (c). Let $\c$ be a decreasing saturated chain from $\hat{0}$ to $X$ in $\NC_G$. By Theorem~\ref{minLabelNBC} part (b), the labels along $\c$ form an NBC set.  Thus, it suffices to show these NBC sets are noncrossing.   Suppose this was not the case and that there were  crossing edges $f_1$ and $f_2$ with $f_1\lhd f_2$ in some  NBC set of $X$ which appears along a saturated chain from $\hat{0}$ to $X$ in $\NC_G$.  We may assume that $X$ is minimal among elements of $\NC_G$ that have crossing edges in one of its NBC sets.

We claim that $f_1$ is the smallest edge of $E(X)$. First note that every spanning NBC set of $X$ contains the smallest edge of $E(X)$. To see why, note that if the smallest edge was a bridge of $X$, it must be in this spanning set.  If it is not a bridge, it is contained in some cycle and so must be in the spanning set since otherwise $X$ would contain a broken circuit.  Thus, the smallest edge of $X$  is in every spanning NBC set of $X$. It follows that the labels along any saturated  chain from $\hat{0}$ to $X$ must contain the smallest edge.   Thus, the labels along every decreasing saturated  chain from $\hat{0}$ to $X$ must end with the smallest edge.  Since $X$ is minimal with respect to having a crossing, the last label must either be $f_1$ or $f_2$. Since $f_1\lhd f_2$, $f_1$ is the smallest edge in $E(X)$.  Now consider the interval $[f_2,X]$ in $\NC_G$.  Then $f_1=\min E(X)\setminus \{f_2\}$ and so by assumption the bond induced on $f_1$ and $f_2$ is a noncrossing bond. Since $f_1$ and $f_2$ are crossing edges, they cannot share an endpoint.  It follows that the bond induced on $f_1$ and $f_2$ is $\{f_1,f_2\}$.  But then the fact that this bond is noncrossing, contradicts the fact that $f_1$ and $f_2$ cross. 

Next, we prove (d).  Let $F=\{e_1,e_2,\dots, e_k\}$ be a spanning noncrossing NBC set of $X$ with  $e_1\rhd e_2\rhd \cdots \rhd e_k$.  Let $H_0=\hat{0}$ and for $1\leq i\leq k$, let $F_i$ be the forest with vertex set $V(G)$ and edge set $\{e_1,e_2,\dots, e_i\}$. Moreover, let  $H_i$ be the bond induced on $F_i$.  We claim that each of these bonds is noncrossing. To see why, note that the partitions associated to $F_j$ and $H_j$ are the same for all $0\leq j\leq k$.  Since each $F_j$ is noncrossing, Proposition~\ref{crossingBondProp} implies that the partition for $F_j$ is noncrossing and so $H_j$ is a noncrossing bond. 

By Theorem~\ref{minLabelNBC}, each spanning NBC set of $X$ appears exactly once as a saturated decreasing chain in $L_G$.  By construction, $H_0\cover H_1\cover \cdots \cover H_k$ is the saturated chain chain from $\hat{0}$ to $X$ which produces the NBC set $F=\{e_1,e_2,\dots, e_k\}$.  The claim now follows since this chain also exists in $\NC_G$.

Finally, for (e),  note that by Proposition~\ref{gradedProp} and part (a), $\rho(H) = n-cc(H)$ for all $H\in \NC_G$. Parts (c) and (d) provide a bijection between the saturated decreasing chains from $\hat{0}$ to $H$ and noncrossing NBC sets.  Then to finish, apply part (b)  and Theorem~\ref{mobFuncEL}.
\end{proof}

After giving this sufficient condition for the minimum labeling to be an EL-labeling, the next obvious question is which graphs satisfy this condition.  We will leave this discussion 
for the next section where we will explore several families of graphs which have this property. In particular, we will show the previous theorem applies to perfectly labeled graphs (Definition~\ref{perLabDef}), upper crossing closed graphs that are tightly closed (Definition~\ref{tightlyClosedDef}), and strongly upper crossed  graphs (Definition~\ref{stronglyUpperCrossedDef}).   The perfectly labeled and strongly upper crossed  graphs include graphs that are not crossing closed (and hence not upper crossing closed).   On the other hand, not every upper crossing closed graph satisfies the conditions of the previous theorem. Because of this, we genuinely require both Theorem~\ref{NCNBCThm} and Theorem~\ref{minLabelNC} to get the results concerning the noncrossing NBC interpretation for the M\"obius function.

We should note that not all graphs $G$ have a shellable $NC_G$. The 5-pointed star in Figure~\ref{5pointStarFig} is one such example, as one can check in SageMath \cite{sagemath}.

 \section{Families of Graphs}\label{famOfGraphsSec}
 In this section, we will consider three families of graphs: perfectly labeled graphs in Subsection~\ref{perfLabSec}, tightly closed graphs in Subsection~\ref{tightClosedSec}, and strongly upper crossed graphs in Subsection~\ref{strongUpCrossSec}.  We present some of the nice properties that the noncrossing bond posets of these families of graphs have such as gradedness, shellability, and combinatorial formulas for M\"obius values.   We finish the section by gathering all the results about the structure of the noncrossing bond poset for the families of graphs studied within this article.  This information can be found in Table~\ref{MetaThmTable}.
\subsection{Perfectly Labeled Graphs }\label{perfLabSec}

In this subsection we present two main results, Theorem~\ref{perfLabMainThm1} and Theorem~\ref{perfLabSuperSol}, concerning perfectly labeled graphs.  We start with the definition of these graphs.

\begin{definition}\label{perLabDef}
Let $G$ be a graph.  We say $G$ is \emph{perfectly labeled} if whenever $ik, jk \in E(G)$ with $i<j<k$, $ij \in E(G)$.  \footnote{We note here that it is common to call a labeling of a graph with this property a \emph{perfect elimination order}.}
\end{definition}

\begin{figure}
\begin{center}
\begin{tikzpicture}[scale=.7]
\coordinate (v1) at (0,1);
\coordinate (v2) at (.87,.5);
\coordinate (v3) at (.87,-.5);
\coordinate (v4) at (0,-1);
\coordinate (v5) at (-.87,-.5);
\coordinate (v6) at (-.87,.5);

\draw[thick] (v6)--(v1) -- (v2) -- (v3)--(v5);
\draw[thick] (v5)--(v6);
\draw[thick] (v1) -- (v4);
\draw[thick] (v1) -- (v3);
\draw[thick] (v1) -- (v5);
\foreach \v in {v1,v2,v3,v4, v5,v6} \fill(\v) circle (.1);
\draw(0,1.4) node{\footnotesize 1};
\draw (.87+.4,.5+.4) node{\footnotesize 2};
\draw(0,-1.4) node{\footnotesize 4};
\draw (.87+.4,-.5-.4) node{\footnotesize 3};
\draw (-.87-.4,-.5-.4) node{\footnotesize 5};
\draw (-.87-.4,.5+.4) node{\footnotesize 6};
\draw(0,-2.2) node {\large $G$};

\begin{scope}[shift = {(8,0)}]

\coordinate (v1) at (0,1);
\coordinate (v2) at (.87,.5);
\coordinate (v3) at (.87,-.5);
\coordinate (v4) at (0,-1);
\coordinate (v5) at (-.87,-.5);
\coordinate (v6) at (-.87,.5);

\draw[thick] (v3)--(v4);
\draw[thick] (v5)--(v4);
\draw[thick] (v1)--(v4);
\draw[thick] (v5)--(v6)--(v2);
\draw[thick] (v2) -- (v3);
\draw[thick] (v4) -- (v6);
\draw[thick] (v4) -- (v2);
\foreach \v in {v1,v2,v3,v4, v5,v6} \fill(\v) circle (.1);
\draw(0,1.4) node{\footnotesize 1};
\draw (.87+.4,.5+.4) node{\footnotesize 2};
\draw(0,-1.4) node{\footnotesize 4};
\draw (.87+.4,-.5-.4) node{\footnotesize 3};
\draw (-.87-.4,-.5-.4) node{\footnotesize 5};
\draw (-.87-.4,.5+.4) node{\footnotesize 6};
\draw(0,-2.2) node {\large $H$};

\end{scope}

\end{tikzpicture}
\end{center}
\caption{Two labeled graphs. $G$ is perfectly labeled whereas $H$ is not.}\label{perfectLabelExample}
\end{figure}
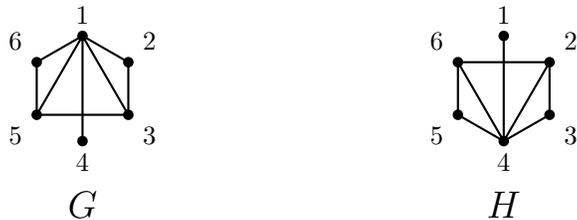

The  graph $G$ in Figure~\ref{perfectLabelExample} is perfectly labeled.  However,  $H$ is  not perfectly labeled, since, for example, $14$ and $24$ are edges, but $12$ is not an edge.   It turns out that there is a classification of graphs which can be perfectly labeled.

\begin{definition}
A graph is \emph{chordal} if for every cycle of length at least 4 there is an edge between two vertices in the cycle which are not adjacent in the cycle.
\end{definition}

The reader may have noticed that  the graphs in Figure~\ref{perfectLabelExample} are both chordal, but only $G$ is perfectly labeled.  It is well-known that a graph can be perfectly labeled if and only if it is chordal (see, for example, the note immediately following Corollary 4.10 of~\cite{s:iha}).  However, as we saw not every labeling of a chordal graph gives rise to a perfectly labeled graph.  The distinction between perfectly labeled and chordal is immaterial to the structure of the bond lattice since the lattice does not depend on the labeling of the vertex set.  However, in the case for the noncrossing bond poset, the structure of the poset can  depend on the labeling of the graph.  Because of this, we focus on perfectly labeled graphs as opposed to just chordal graphs.   This is also why we use the term ``perfectly labeled" as opposed to saying ``$G$ has a perfect elimination order" which is more common in the literature.

As we will see  throughout this subsection, increasing spanning trees and forests play an important role of the combinatorics of the noncrossing bond posets of perfectly labeled graphs.  

\begin{definition}
Let $T$ be a tree with vertices which are distinct integers.  Let $r$ be the smallest vertex of $T$. We say $T$ is an \emph{increasing tree} if the vertices along any path from $r$ to any other vertex form an increasing sequence.   We say a spanning subgraph of a graph $G$ is an \emph{increasing spanning forest} of $G$ if each connected component is an increasing tree. 
\end{definition}

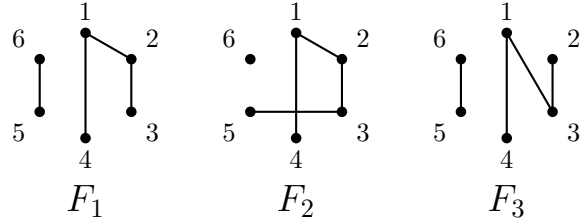
\begin{figure}
\begin{center}
\begin{tikzpicture}[scale=.7]
\coordinate (v1) at (0,1);
\coordinate (v2) at (.87,.5);
\coordinate (v3) at (.87,-.5);
\coordinate (v4) at (0,-1);
\coordinate (v5) at (-.87,-.5);
\coordinate (v6) at (-.87,.5);

\draw[thick] (v4)--(v1)--(v2)--(v3);
\draw[thick] (v5)--(v6);
\foreach \v in {v1,v2,v3,v4, v5,v6} \fill(\v) circle (.1);
\draw(0,1.4) node{\footnotesize 1};
\draw (.87+.4,.5+.4) node{\footnotesize 2};
\draw(0,-1.4) node{\footnotesize 4};
\draw (.87+.4,-.5-.4) node{\footnotesize 3};
\draw (-.87-.4,-.5-.4) node{\footnotesize 5};
\draw (-.87-.4,.5+.4) node{\footnotesize 6};
\draw(0,-2.2) node {\large $F_1$};

\begin{scope}[shift = {(4,0)}]

\coordinate (v1) at (0,1);
\coordinate (v2) at (.87,.5);
\coordinate (v3) at (.87,-.5);
\coordinate (v4) at (0,-1);
\coordinate (v5) at (-.87,-.5);
\coordinate (v6) at (-.87,.5);

\draw[thick] (v1)--(v2)--(v3)--(v5);
\draw[thick] (v1)--(v4);
\foreach \v in {v1,v2,v3,v4, v5,v6} \fill(\v) circle (.1);
\draw(0,1.4) node{\footnotesize 1};
\draw (.87+.4,.5+.4) node{\footnotesize 2};
\draw(0,-1.4) node{\footnotesize 4};
\draw (.87+.4,-.5-.4) node{\footnotesize 3};
\draw (-.87-.4,-.5-.4) node{\footnotesize 5};
\draw (-.87-.4,.5+.4) node{\footnotesize 6};
\draw(0,-2.2) node {\large $F_2$};

\end{scope}

\begin{scope}[shift = {(8,0)}]

\coordinate (v1) at (0,1);
\coordinate (v2) at (.87,.5);
\coordinate (v3) at (.87,-.5);
\coordinate (v4) at (0,-1);
\coordinate (v5) at (-.87,-.5);
\coordinate (v6) at (-.87,.5);

\draw[thick] (v4)--(v1)--(v3)--(v2);
\draw[thick]  (v5)--(v6);
\foreach \v in {v1,v2,v3,v4, v5,v6} \fill(\v) circle (.1);
\draw(0,1.4) node{\footnotesize 1};
\draw (.87+.4,.5+.4) node{\footnotesize 2};
\draw(0,-1.4) node{\footnotesize 4};
\draw (.87+.4,-.5-.4) node{\footnotesize 3};
\draw (-.87-.4,-.5-.4) node{\footnotesize 5};
\draw (-.87-.4,.5+.4) node{\footnotesize 6};
\draw(0,-2.2) node {\large $F_3$};

\end{scope}

\end{tikzpicture}
\end{center}
\caption{Several Forests  }\label{isfExample}
\end{figure}

\noindent The forests $F_1$ and $F_2$ in Figure~\ref{isfExample}  are increasing spanning forests of the graph $G$ in Figure~\ref{perfectLabelExample}, whereas  $F_3$ is not since the path from $1$ to $2$ is not increasing.

We now present a series of lemmas aimed at showing that perfectly labeled graphs satisfy the conditions of Theorem~\ref{minLabelNC}. This will allow us to show that the noncrossing bond poset of a perfectly labeled graph is graded, shellable, and has a combinatorial interpretation for the M\"obius function and characteristic polynomial.

\begin{lem}\label{incPath}
Let $G$ be a connected perfectly labeled graph. Let $r$ be the smallest vertex in $G$ and let $v$ be any other vertex of $G$.  There exists an increasing path from $r$ to $v$ in $G$.  That is, there is a path $ru_1u_2\dots u_k v$ where $r<u_1<u_2<\cdots <u_k<v$.
\end{lem}

\begin{proof}
Suppose this was not the case.  Let $P$ be a shortest path from $r$ to $v$. By assumption, $P$ is not an increasing path. Thus it must contain contain a sequence of vertices $ikj$ where $i<j<k$.  But then since $G$ is perfectly labeled, $ij\in E(G)$, so we may replace the edges $ik, jk$ with $ij$ and get a shorter path from $r$ to $v$.  This contradicts the minimality of $P$.  
\end{proof}

\begin{lem}\label{perfLabLem}
Let $G$ be a perfectly labeled graph. Suppose that $H\leq H'$ in $\NC_G$.  Moreover, suppose that $B_1,B_2,\dots, B_k$ where $\min B_1<\min B_2<\cdots < \min B_k$ are the connected components of $H$ that are merged together to get $H'$.  Then merging $B_1$ and $B_2$ in $H$ creates a noncrossing bond of $G$.
\end{lem}
\begin{proof}
Let $3\leq i \leq k$.  If merging $B_1$ and $B_2$ crossed with some  $B_i$, then there exists $a,c \in B_1 \cup B_2$ and $b,d \in B_i$ or $a,c \in B_i$ and $b,d \in B_1\cup B_2$ with $a<b<c<d$.  If $a,c \in B_1 \cup B_2$ and $b,d \in B_i$, then $\min B_1< \min B_2 <b < c<d$ which implies either $B_1$ and $B_i$ cross or $B_2$ and $B_i$ cross (depending on if $c\in B_1$ or $c\in B_2$).  Neither is possible since $H$ is noncrossing.  A similar argument shows that it is not possible that there exists $a,c \in B_i$ and $b,d \in B_1\cup B_2$ with $a<b<c<d$. Hence merging $B_1$ and $B_2$ does not cause a crossing with $B_3,B_4,\dots, B_k$. Moreover, merging $B_1$ and $B_2$ cannot create a crossing with any of other connected components of $G$ since that would mean that $H'$ was crossing.  

Thus it suffices to show that merging $B_1$ and $B_2$ in $H$ forms a bond of $G$. 
It is not hard to see that any induced subgraph of a perfectly labeled graph is perfectly labeled.  It follows that $H'$ is perfectly labeled. By Lemma~\ref{incPath} there is an increasing path in $H'$  from $\min B_1$ to $\min B_2$.  Except for $\min B_2$, this path must only contain vertices from $B_1$ since otherwise it would not be increasing. So there is an edge, $e$, between a vertex in $B_1$ and a vertex in $B_2$.  Then the bond induced on on $E(H)\cup \{e\}$ is exactly the spanning subgraph obtained by merging $B_1$ and $B_2$ in $H$. The result now follows.
\end{proof}

In the following lemma we will order the edges \emph{colexicographically}.  That is, we say $ab \lhd a'b'$ if and only if $b<b'$ or $b=b'$ and $a < a'$.  Note that colexicographic order is just lexicographic order reading right to left instead of left to right.

\begin{lem}\label{perLabMinEdgeLem}
Let $G$ be a perfectly labeled graph with the edges colexicographically.  If $H<H'$ and $e=\min E(H')\setminus E(H)$, then the bond induced on $E(H)\cup \{e\}$ is noncrossing.
\end{lem}
\begin{proof}
Let $B_1,B_2,\dots, B_k$ be the blocks which merge together to get $H'$ where  $\min B_1<\min B_2<\cdots < \min B_k$.  Since $G$ is perfectly labeled and $H'$ is an induced subgraph, it is perfectly labeled. By Lemma~\ref{incPath}, in the connected component containing $B_1$ and $B_2$ there is an increasing path from $\min B_1$ to $\min B_2$.  The last edge on this path must be of the form $a\min B_2$ where $a\in B_1$.  Let $a'$ be the smallest vertex in $B_1$ for which there is an edge $a'\min B_2$.  Since we are ordering edges colexicographically, $a'\min B_2=\min E(H')\setminus E(H)$.    The bond induced on $E(H)\cup \{a'\min B_2\}$ is the bond obtained by merging $B_1$ and $B_2$.  By Lemma~\ref{perfLabLem} this bond is noncrossing.
\end{proof}

As with noncrossing NBC sets, we say an increasing spanning forest is \emph{noncrossing} if none of the edges cross. For example, the forests $F_1$ and $F_2$ in Figure~\ref{isfExample} are increasing spanning forests of $G$ in Figure~\ref{perfectLabelExample}, but only $F_1$ is noncrossing. It turns out that when we use the colexicographic order, there are the same number of noncrossing increasing spanning forests and noncrossing NBC sets. We use the notation $\ncisf_k(G)$ for the number of noncrossing increasing spanning forest of $G$ with $k$ edges.  In~\cite[Theorem 2.4]{hms:isf}, it was shown that when $G$ is perfectly labeled and we order the edges lexicographically, the NBC sets are exactly the increasing spanning forests.  The proof given there can easily be modified to allow for the case when the edges are ordered colexicographically.  Thus we have the following.

\begin{lem}\label{equalityNBCisfLem}
Let $G$ be a perfectly labeled graph with the edges ordered colexicographically.  Then for all $k\geq 0$,
$$
\ncisf_k(G) =\ncnbc_k(G).
$$
\end{lem}

We now present the first main theorem of this section.  Applying Lemma~\ref{perLabMinEdgeLem}, Lemma~\ref{equalityNBCisfLem}, and Theorem~\ref{minLabelNC} we get the following.

\begin{thm}\label{perfLabMainThm1}
Let $G$ be a perfectly labeled graph on $[n]$ such that $\NC_G$ has a $\hat{1}$.  Then we have the following.
\begin{enumerate}
    \item [(a)] $\NC_G$ is graded and for $H\in \NC_G$, the rank of $H$ is given by $\rho(H) = n -cc(H)$.
    \item[(b)]The minimum labeling with respect to the colexicographic ordering on $E(G)$ is an EL-labeling of $\NC_G$ and so $\NC_G$ is shellable.
    \item [(c)] For $H\in\NC_G$,
    $$
\mu(H) = (-1)^{n-cc(H) } \ncnbc_{n-cc(H)}(H) =  (-1)^{n-cc(H) } \ncisf_{n-cc(H)}(H) 
$$
and
$$
\chi(\NC_G, t) =\sum_{k\geq 0} (-1)^k\ncnbc_k(G) t^{\rho(\NC_G) -k} =\sum_{k\geq 0} (-1)^k\ncisf(G) t^{\rho(\NC_G) -k},
$$
where the NBC sets are with respect to the colexicographic ordering on $E(G)$.        
\end{enumerate}
\end{thm} 

We note that in~\cite{bs:mfl}, Blass and Sagan used NBB sets to show that the  M\"obius function of  the noncrossing partition lattice counts noncrossing increasing trees and hence is a Catalan number.  The previous theorem generalizes this result since the noncrossing partition lattice is the noncrossing bond poset of the complete graph which is perfectly labeled. 

\begin{figure}
\begin{center}
\begin{tikzpicture}[scale=.7]
\coordinate (v1) at (0,1);
\coordinate (v2) at (1,0);
\coordinate (v3) at (0,-1);
\coordinate (v4) at (-1,0);
\draw[thick] (v3)--(v1) -- (v2) -- (v4);
\foreach \v in {v1,v2,v3,v4} \fill(\v) circle (.1);
\draw(0,1.4) node{\footnotesize 1};
\draw(0,-1.4) node{\footnotesize 3};
\draw(1.4,0) node{\footnotesize 2};
\draw(-1.4,0) node{\footnotesize 4};
\draw(0,-3.2) node {\large $G$};

\begin{scope}[shift = {(8,-2)}]

\node (1/2/3/4) at (0,0)  {$1/2/3/4$};

\node (12/3/4) at (0,2)  {$12/3/4$};
\node (13/2/4) at (-5,2)  {$13/2/4$};
\node (1/24/3) at (5,2)  {$1/24/3$};

\draw[thick] (1/2/3/4)--(12/3/4);
\draw[thick] (1/2/3/4)--(13/2/4);
\draw[thick] (1/2/3/4)--(1/24/3);

\node (123/4) at (-5,4)  {$123/4$};
\draw[thick] (12/3/4)--(123/4)--(13/2/4);
\node (124/3) at (5,4)  {$124/3$};
\draw[thick] (12/3/4)--(124/3)--(1/24/3);

\node (1234) at (0,6)  {$1234$};
\draw[thick] (123/4)--(1234)--(124/3);
\draw(0,-1.2) node {\large $\NC_G$};
\tikzstyle{every node}=[blue, scale=1.1]
\node at (-2.5, 1.25) {$2$};
\node at (0.2, 1.2) {$1$};
\node at (2.5, 1.25) {$3$};
\node at (-2.5, 3.25) {$2$};
\node at (2.5, 3.25) {$3$};
\node at (5.5, 3.25) {$1$};
\node at (-5.5, 3.25) {$1$};
\node at (-2.5, 5.25) {$3$};
\node at (2.5, 5.25) {$2$};
\end{scope}
\end{tikzpicture}
\end{center}
\caption{A graph and its noncrossing bond poset labeled by the max-min edge labeling. }\label{SnLabelFig}
\end{figure}
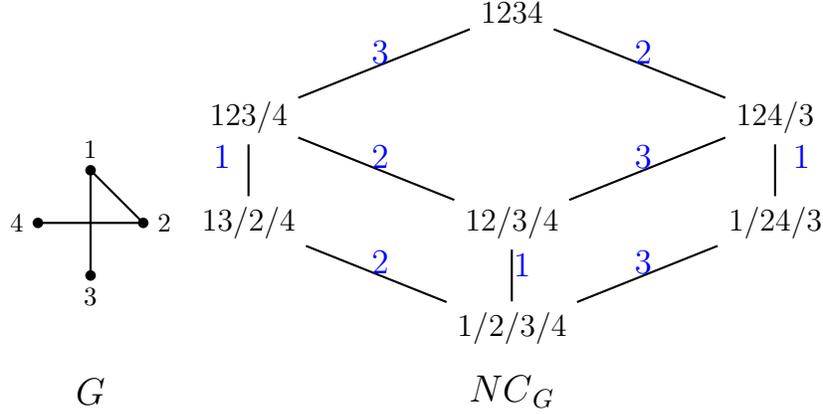

 Let us now turn our attention to other the main result of this subsection.  In~\cite{s:sl} Stanley introduced the notion of a supersolvable lattice.  A lattice $L$ is called \emph{supersolvable} if there exists a maximal chain with the property that it along with any other chain in $L$ generates (by taking joins and meets of all elements in the chains) a distributive lattice.   It is well-known that a graph is chordal if and only if its bond lattice is supersolvable (see, for example,  Corollary 4.10 and the note that follows it in~\cite{s:iha}).  As we will see in Theorem~\ref{perfLabSuperSol}, there is an analogue of this result for noncrossing bond posets.  However,  we should point out that when we pass to the noncrossing case things become a bit more complicated. First, as we noted earlier, the noncrossing bond poset depends on the labeling of the vertices.  Because of this there are chordal graphs which are not perfectly labeled and such that the noncrossing bond poset is not supersolvable.  Second, we no longer have an if and only if statement as there exist noncrossing bond posets that are supersolvable lattices, but do not come from chordal graphs. Finally, not every perfectly labeled graph is crossing closed and so the noncrossing bond poset of a perfectly labeled graph may  not even be a lattice.

In order to study supersolvability of the noncrossing bond poset we will consider a special type of EL-labeling of graded posets. Let $\lambda$ be an EL-labeling of a poset of rank $n$.  We say $\lambda$ is  an \emph{$S_n$ EL-labeling} if every maximal chain of $P$ is labeled by a permutation of $[n]$ (with natural ordering on $[n]$).  The edge labeling in Figure~\ref{SnLabelFig} is an example of an $S_n$ EL-labeling.  We note  that the condition that the unique maximal chain in each interval is lexicographically first is automatically implied if the maximal chains are labeled by permutations and thus, we  only need to check  that each interval has a unique    increasing maximal chain. In addition to the properties that EL-labeled posets possess, there are special properties that a lattice with an $S_n$ EL-labeling possesses.   In particular, McNamara~\cite{m:els} showed that if $L$ is a graded lattice then $L$ is supersolvable if and only if it has an $S_n$ EL-labeling.  We will show that if $G$ is perfectly labeled  and connected, then it has an $S_n$ EL-labeling.   To do this, we use a labeling introduced by Bj\"orner and Edelman.

\begin{definition}[Bj\"orner-Edelman~\cite{b:scmpos}]
Let $G$ be a graph on $[n]$ such that $\NC_G$ is graded.  The \emph{max-min} edge labeling is defined by
$$
\lambda(H\cover H') = \max\{\min B, \min B'\} -1
$$
where $B$ and $B'$ are the blocks merged  when going from $H$ to $H'$.
\end{definition}

See the poset in Figure~\ref{SnLabelFig} for an example of this labeling.  Bj\"orner and Edelman~\cite{b:scmpos} showed that the max-min edge labeling gives an EL-labeling of the noncrossing partition lattice (which is also the noncrossing bond poset of the complete graph).  As McNamara points out in~\cite{m:els}, this labeling is in fact an $S_n$ EL-labeling.  It turns out that among connected graphs, perfectly labeled graphs are exactly the graphs where the max-min edge labeling is an $S_n$ EL-labeling of the noncrossing bond poset.   Note that in the hypothesis of the following proposition, we assume $G$ has $n+1$ vertices and is connected.  This guarantees that $\NC_G$ has  rank $n$.

\begin{prop}\label{minMaxPerfectLabel}
Let $G$ be a  connected graph on $[n+1]$ such that $\NC_G$ is graded.  The max-min edge labeling is an $S_n$ EL-labeling of $\NC_G$ if and only if $G$ is perfectly labeled.
\end{prop}

\begin{proof}
  
Given that $G$ has $n+1$ vertices, $\NC_G$ is a subposet of $\NC_{K_{n+1}}$. Since $\NC_{K_{n+1}}$ is the noncrossing partition lattice on $[n+1]$, we have that $\NC_G$ is a subposet of the noncrossing partition lattice.  It follows from the assumption that $\NC_G$ is graded with the same rank function as the noncrossing partition lattice (see Proposition~\ref{gradedProp}), we have that the set of maximal chains in $\NC_G$ is a subset of the maximal chains of the noncrossing partition lattice. Since the cover relations in $\NC_G$ are the same as that in the noncrossing partition lattice, we also have that the label sequences that appear along the maximal chains in $\NC_G$ are the same as those that appear in the noncrossing partition lattice.  Since  it is known that the max-min edge labeling is an $S_n$ EL-labeling of the noncrossing partition lattice, to finish the  proof we can show that  each interval of $\NC_G$ has an   increasing  maximal chain if and only if $G$ is perfectly labeled.

Suppose that $G$ is not perfectly labeled.  Then there exists edges $ik, jk$ such that $i<j<k$  and $ij \notin E(G)$.  Let $H$ be the bond of $G$ where $i,j,k$ are in the same connected component and every other connected component is trivial.  Consider the interval $[\hat{0}, H]$.  Since $ij\notin E(G)$, this interval has two maximal chains both labeled by $k-1, j-1$.  Thus, the interval has no increasing chain.

Next, suppose that $G$ is perfectly labeled.   Suppose that $[X,Y]$ is an interval in $\NC_G$ and suppose that $B_1,B_2,\dots, B_k$ are the connected components of $X$ that will merge together to get $Y$. Moreover, assume that $\min B_1<\min B_2<\cdots < \min B_k$.  It is not hard to see that if there is an  increasing maximal chain  in $[X,Y]$, the first step must be to merge $B_1$ and $B_2$.   Let $Z$ be the bond obtained by merging $B_1$ and $B_2$ in $X$. We can apply Lemma~\ref{perfLabLem} to see that $Z \in \NC_G$.   Now we can use induction to prove that $[Z,Y]$ has an   increasing  maximal chain which can be concatenated with the label from $X$ to $Z$ to give an   increasing maximal chain in $[X,Y]$.
\end{proof}

In Proposition~\ref{minMaxPerfectLabel}, we assumed that $G$ is connected. By Lemma \ref{lem:product}, if $G$ is not connected and its connected components do not cross, $\NC_G$ is the product of smaller noncrossing bond posets, one for each connected component.  If $G$ is perfectly labeled, each connected component of $G$ must be perfectly labeled.  Thus, if $G$ perfectly labeled, the noncrossing bond poset of each of its connected components has an $S_n$ EL-labeling.  Using McNamara's~\cite{m:els} result about $S_n$ EL-labelings and supersolvability implies each is noncrossing bond poset is supersolvable.  Moreover, if $G$ is crossing closed its connected components cannot cross.  Putting this altogether and using the fact that the product of supersolvable lattices is a supersolvable lattice, we get the following.

\begin{thm}\label{perfLabSuperSol}
Let $G$ be a perfectly labeled graph.  If $G$ is crossing closed, then $\NC_G$ is a supersolvable lattice.
\end{thm}

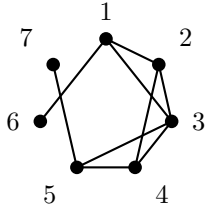
\begin{figure}
\begin{center}
\begin{tikzpicture}[scale=.9]
\coordinate (v1) at (0,1);
\coordinate (v2) at (.78,.62);
\coordinate (v3) at (.97,-0.22);
\coordinate (v4) at (0.43,-0.9);
\coordinate (v5) at (-0.43,-0.9) ;
\coordinate (v6) at (-.97,-0.22);
\coordinate (v7) at (-.78,.62);

\draw[thick] (v1) -- (v2)--(v3)--(v1)--(v6);
\draw[thick] (v3)--(v4)--(v5)--(v3);
\draw[thick] (v5)--(v7);
\draw[thick] (v2)--(v4);

\foreach \v in {v1,v2,v3,v4, v5,v6, v7} \fill(\v) circle (.1);

\draw(0,1.4) node{\footnotesize 1};
\draw(.78+.4,.62+.4) node{\footnotesize 2};
\draw(.97+.4,-0.22) node{\footnotesize 3};
\draw (0.43+.4,-0.9-.4) node{\footnotesize 4};
\draw (-0.43-.4,-0.9-.4) node{\footnotesize 5};
\draw (-.97-.4,-0.22) node{\footnotesize 6};
\draw (-.78-.4,.62+.4) node{\footnotesize 7};

\end{tikzpicture}
\end{center}
\caption{A perfectly labeled graph which is not crossing closed\label{peoNotCCFig} }
\end{figure}

We mention here that not every  perfectly labeled graph is crossing closed (hence the necessity of the crossing closed hypothesis in Theorem~\ref{perfLabSuperSol}).  The graph  in Figure~\ref{peoNotCCFig} is perfectly labeled, but not crossing closed.  This is because there are two minimal induced connected components containing $16$ and $57$, namely  the one containing the vertices $1,3,5,6,7$ and the one containing $1,2,4,5, 6,7$.  Nevertheless, when $G$ is crossing closed and perfectly labeled, it is upper crossing closed as we see next.

\begin{prop}
Let $G$ be a perfectly labeled graph which is crossing closed.  Then $G$ is upper crossing closed with respect to the colexicographic and lexicographic order.
\end{prop}
\begin{proof}
Suppose that $ac$ and $bd$ cross with $a<b<c<d$.  Then in $J(ac,bd)$, there is a path from $a$ to $b$.  Let $P:av_1v_2\dots v_kb$ be a path from $a$ and $b$ which is minimal with respect to length.  If $P$ is not increasing, then there is a an index $i$ with $v_{i-1},v_{i+1}<v_{i}$.  But then since $G$ is perfectly labeled, there is an edge $v_{i-1}v_{i+1}$ contradicting the minimality of $P$.  So $P$ must be increasing.  Then $av_1$ is an edge in $E(G)$ and is smaller in colexicographic and lexicographic order than $ac$ and $bd$.  It follows that $G$ is upper crossing closed with respect to colexicographic and lexicographic order.
\end{proof}

The reader may be wondering if $\NC_G$ being supersolvable implies that $G$ is chordal since this is the case for the bond lattice of a graph.  The graph in Figure~\ref{twistedC4Fig} shows this is not true.  It is a 4-cycle and thus is not chordal.  Nevertheless, its noncrossing bond poset is a supersolvable lattice.

\subsection{Tightly Closed Graphs}\label{tightClosedSec}
As we saw in Section~\ref{structureSec}, crossing closed graphs need not have graded noncrossing bond posets.  However, it turns out that if we restrict what $J(e,f)$ can look like, we can guarantee the noncrossing bond poset is graded.  Moreover,  if we further assume the graph is upper crossing closed, we obtain more properties of the poset.  We explore these ideas next.

\begin{definition}\label{tightlyClosedDef}
Let $G$ be a graph. We say $G$ is \emph{tightly closed} if it is crossing closed and for all edges $e$ and $f$ that cross, $J(e,f)$ is a subgraph of $K_4$.
\end{definition}

  The complete bipartite graphs are a family  of tightly closed graphs.  To see why, note that if two edges cross in a complete bipartite graph, they must connect the two parts of the graph and so must lie on a (twisted) 4-cycle.  The 5-pointed star depicted in Figure~\ref{5pointStarFig} gives a different an example of a tightly closed graph.  Note that since the 5-pointed star is a cycle, it is 2-connected.   It turns out that any 2-connected crossing closed graph is tightly closed.

\begin{prop}
If $G$ is 2-connected and crossing closed, then $G$ is tightly closed.
\end{prop}
\begin{proof}
By Lemma~\ref{mefLemma}, if $G$ was not tightly closed, $G$ would have cut vertices.    This is impossible as $G$ is 2-connected.
\end{proof}

\begin{figure}
\begin{center}
\begin{tikzpicture}[scale=.9]
\coordinate (v1) at (0,1);
\coordinate (v2) at (.87,.5);
\coordinate (v3) at (.87,-.5);
\coordinate (v4) at (0,-1);
\coordinate (v5) at (-.87,-.5);
\coordinate (v6) at (-.87,.5);

\draw[thick] (v5)-- (v4)--(v1)--(v2)--(v3)--(v6)--(v5);

\foreach \v in {v1,v2,v3,v4, v5,v6} \fill(\v) circle (.1);
\draw(0,1.4) node{\footnotesize 1};
\draw (.87+.4,.5+.4) node{\footnotesize 2};
\draw(0,-1.4) node{\footnotesize 4};
\draw (.87+.4,-.5-.4) node{\footnotesize 3};
\draw (-.87-.4,-.5-.4) node{\footnotesize 5};
\draw (-.87-.4,.5+.4) node{\footnotesize 6};
\draw(0,-2.2) node {\large $G$};

\end{tikzpicture}
\end{center}
\caption{A 2-connected graph that is not crossing closed \label{twistedC6} }
\end{figure}
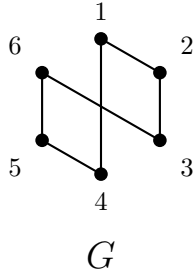

Before we move on, we wish to mention that 2-connected does not imply crossing closed.  For example, the 6-cycle in Figure~\ref{twistedC6} is 2-connected, but not crossing closed.

\begin{thm}\label{tightlyClosedGraded}
If $G$ be a tightly closed graph, then $\NC_G$ is graded.
\end{thm}
\begin{proof}
Let $H\cover H'$.   For each edge $f\in E(H')\setminus E(H)$, we will count crossings of $f$ with the edges in $E(H)$.  We will call such a crossing bad if the edges in the crossing are in different components of the graph with edge set $E(H) \cup \{f\}$ and vertex set $V(G)$. 

We claim that there is an edge in $E(H')\setminus E(H)$ with no bad crossings. To see why suppose that this was not the case and let $e=ac$ be an edge of $E(H')\setminus E(H)$ with a minimum number of bad crossings. By assumption, $e$ has at least one bad crossing, say with edge $e'=bd$ in $H$ where $a<b<c<d$.    Since $e$ and $e'$ are in $H'$ and are crossing and $H'$ is noncrossing, they must lie in the same connected induced component of $H'$.  Thus $J(e,e')$ is a subgraph of this component.  Since $G$ is tightly closed, $J(e,e')$ is a subgraph of $K_4$ and so one of the edges $ab, bc, cd, ad$ must be present in $H'$.  Without loss of generality we may assume that $ab$ is present.  Note that $ab$ is not in $H$ as $e$ and $e'$ is assumed to be a bad crossing and $e$ and $e'$ thus lie in different components of the graph with edge set $E(H) \cup \{e\}$.

Since there are no edges in $E(H')\setminus E(H)$ with no bad crossings, $ab$ must have a bad crossing with some edge, say $vw$.  We claim that $vw$ must cross $e$. If this was not the case, then $a<v<b<w\leq c<d$.  This implies that $vw$ crosses $e' = bd$ in $H$ and since $H$ is noncrossing, $vw$ and $e' =  bd$ must be in the same component. But then $vw$ and $ab$ do not form a bad crossing. Thus, $vw$ crosses $e$ and so any edge that crosses $ab$ to form a bad crossing will form a bad crossing with $e$.  This means that $e$ has strictly more bad crossings than $ab$ which is impossible as $e$ was chosen to be minimal.  Thus, there is an edge of $E(H')\setminus E(H)$ which has no bad crossings.

Let $e$ be an element of $ E(H')\setminus E(H)$ which has no bad crossings and let $G'$ be the graph on $V(G)$ with edge set $E(H) \cup \{e\}$.  Since $e$ has no bad crossings, the partition associated to $G'$ is noncrossing.  Now let $H''$  be  the bond induced on $E(H) \cup \{e\}$.  Then $G'$ and $H''$ have the same connected components and so correspond to the same partition.  It follows that $H''$ is a noncrossing partition.  Since $H<H''\leq H'$ and $H\cover H'$, $H'=H''$.  Moreover, by construction there are exactly two components that merge together from $H$ to $H'$ and so  by Proposition~\ref{gradedProp}, $\NC_G$ is graded.
\end{proof}

The notion of tightly closed may seem artificial at first glance.  However, there is an order-theoretic way to define tightly closed that is much like the notion of semimodularity.  Recall that a lattice $L$ is (upper) semimodular, if it is graded and for all $x,y\in L$, $\rho(x\vee y) +\rho(x\wedge y) \leq \rho(x)+\rho(y)$.  In the case that $a_1$ and $a_2$ are distinct atoms of $L$, semimodularity implies that $\rho(a_1\vee a_2)=2$.  Tightly closed graphs can be defined by slightly relaxing this idea.

\begin{thm}
Let $G$ be a graph which is crossing closed.  $G$ is tightly closed if and only if for all distinct atoms $a_1,a_2\in \NC_G$, $\rho(a_1\vee a_2)=2$ or $\rho(a_1\vee a_2)=3$.
 \end{thm}
 
 \begin{proof}
First note that the atoms of $\NC_G$ are the edges of $G$.  Let $e,f\in E(G)$. If $e$ and $f$ do not cross, $\rho(e\vee f)=2$.  If $e$ and $f$ do cross, then by Proposition~\ref{latticeProps} part (c),  $G$ being crossing closed implies that $e\vee f$ is the bond with a unique nontrivial connected component $J(e,f)$.  Thus, the fact that for any $H\in \NC_G$  $\rho(H) = |V(G)|-cc(H)$ implies that  $\rho(a_1\vee a_2)=3$ if and only if $J(e,f)$ is a connected graph on 4 vertices.  The result now follows.
 \end{proof}

While the noncrossing bond poset of a tightly closed graph is graded, other nice order-theoretic properties do not always hold.  For example, they need not be shellable or have a  noncrossing NBC interpretation.  This can be seen by noting that the 5-pointed star is tightly closed, but it is not shellable and also does not have the NCNBC interpretation for its M\"obius function.  However, if we make the further assumption that the graph is upper crossing closed, we get many nice properties.

\begin{prop}\label{tightlyAndUpperCrossingClosed}
Let $G$ be an upper crossing closed graph which is tightly closed.  Then the conclusions of Theorem~\ref{minLabelNC} hold.  In particular, $\NC_G$ is graded, shellable, and the M\"obius function and characteristic polynomial have  a combinatorial interpretation in terms of NCNBC sets.
\end{prop} 

\begin{proof}
  First note that since $G$ is crossing closed, Proposition~\ref{latticeProps} part (e) implies that $\NC_G$ has a $\hat{1}$. Now suppose that $H<H'$. Let $e=\min E(H')\setminus E(H)$. Adopting the terminology from the proof of Theorem~\ref{tightlyClosedGraded}, we will show that $e$ has no bad crossings in $H\cup \{e\}$.  Suppose this was not the case.  Then there is an edge  $f$ not in the same connected component of $e$ in $E(H)\cup \{e\}$ which crosses $e$.  Since $H'$ is noncrossing and $G$ is tightly closed, there is an edge, $h$, connecting an endpoint of $e$ and an endpoint of $f$.  Moreover, we may assume that $h$ precedes $e$ in $\unlhd$ as $G$ is upper crossing closed.  Since $e=\min E(H')\setminus E(H)$, this would imply that $h\in E(H)$. So $e$ and $f$ do not form a bad crossing, a contradiction.  Since $e$ has no bad crossings, the bond induced on $E(H)\cup \{e\}$ is noncrossing.  Applying Theorem~\ref{minLabelNC} now completes the proof.
\end{proof}

 Before we finish this subsection, let us give an example of  a tightly closed upper crossing closed graph.  Let $K_{even,odd}^n$ be the complete bipartite graph on $[n]$ whose parts are the even and odd numbers. As mentioned at the beginning of this subsection, any complete bipartite graph is tightly closed, so  $K_{even,odd}^n$ is tightly closed. 
Recall that our graphs lie on a circle with evenly spaced vertices.  Given vertices $x<y$ of $K_{even,odd}^n$, we let $\dist(x,y) = \min(y-x-1, n-y+x-1)$ be the minimum number of vertices between $x$ and $y$.  Define a partial order on $E(K_{even,odd}^n)$  by declaring $ij$ is less than  $i'j'$ if and only if $\dist(i,j) < \dist(i',j')$ and let $\unlhd$ be any linear extension of this order. It is not hard to show $K_{even,odd}^n$ is upper crossing closed with respect to this ordering.

\subsection{Strongly Upper Crossed Graphs}\label{strongUpCrossSec}

In this subsection we consider another family of (not necessarily  crossing closed) graphs and show their noncrossing bond posets are graded and shellable.  We also show the noncrossing NBC set interpretation for the M\"obius function and characteristic polynomial still hold in this setting.

\begin{definition}\label{stronglyUpperCrossedDef}
Let $G$ be a graph with a total ordering, $\unlhd$, on the edge set of $G$.  We say that a graph $G$ is \emph{strongly upper crossed} with respect to $\unlhd$ if whenever  $ac, bd$ are crossing edges, there is at least one minimal induced connected component of $G$ containing $ac$ and $bd$ and every edge in each minimal induced connected component of $G$ containing $ac$ and $bd$ precedes $ac$ and $bd$ in the ordering $\unlhd$.
\end{definition}

One may think of strongly upper crossed graphs as a relaxing  of the crossing closed condition, but at the cost of requiring a  stronger  condition  on the ordering of edges as compared to that given for upper crossing closed graphs.  As an example of a strongly upper crossed graph, consider  the graph $G$ in Figure~\ref{graphFig}.  If we order the edges so that  $14$ and $35$ are the largest, $G$ is strongly upper crossed  with respect to this order.  More generally, any connected graph with a single crossing is strongly upper crossed  with ordering where the crossing edges are ordered so that they are the largest.  Since the graph in Figure~\ref{graphFig} is not chordal nor crossing closed, the family of  strongly upper crossed graphs is distinct from the families presented in the previous two subsections.

We should  note that not all upper crossing closed graphs are  strongly upper crossed.  For example, $K_5$ is upper crossing closed with respect to the lexicographic ordering, but there is no ordering which makes it strongly upper crossed.  To see why, note that since $14$ and $25$ cross and $24 \in J(14,25)$, $24$ must be smaller than $14$ and $25$.  However, $24$ and $35$ cross and $25\in J(24,35)$ so $25$ must be smaller than $24$ which is impossible.

The following lemma will allow us to apply Theorem~\ref{minLabelNC} to strongly upper crossed graphs.

\begin{lem}\label{scLem}
Let $G$ be a strongly upper crossed  graph.  If $H< H'$ and $e=\min E(H')\setminus E(H)$, then the bond induced by $E(H)\cup \{e\}$ is noncrossing.
\end{lem}
\begin{proof}
Let $H''$ be the bond induced  by $E(H)\cup \{e\}$ and suppose that  $H''$ is a crossing bond.    Let $B_1$ and $B_2$  be the blocks that are merged when moving from $H$ to $H''$.
Since $H''$ is crossing there is some $f\in E(H)$ which crosses an edge between $B_1$ and $B_2$ and is not in $B_1$ and $B_2$. Thus, $f$ separates $B_1$ and $B_2$.  Since $e$ connects $B_1$ and $B_2$, it too must cross $f$.  Since $H'$ is noncrossing and $G$ is strongly upper crossed, there is a  minimal induced connected component containing $e$ and $f$ in $H'$.  Since $G$ is strongly upper crossed, all the edges in this  minimal induced connected component are smaller than $e$ and $f$.  Not all these edges can be in $H$ since this would imply $f$ did not cause $H''$ to be crossing.  But this is impossible since since $e$ was the smallest edge.
\end{proof}

Using the previous lemma and the fact that noncrossing bond posets of strongly upper crossed graphs have a $\hat{1}$, we get the main theorem of this subsection.

\begin{thm}\label{stronglyUpperCrossedMainThm}
Let $G$ be a strongly upper crossed graph.  Then the conclusions of Theorem~\ref{minLabelNC} hold.  In particular, $\NC_G$ is graded, shellable, and the M\"obius function and characteristic polynomial have  a combinatorial interpretation in terms of NCNBC sets.

\end{thm} 

\subsection{A Summary of Results on Families of Graphs}
\begin{table}[ht]
  \centering
     \begin{tabular}{|c|c|c|c|c|}
    \hline
          & \textbf{Graded} & \textbf{Lattice} & \textbf{\begin{tabular}{@{}c@{}} NCNBC \\ interpretation \end{tabular}} & \textbf{Shellable} \\
    \hline\hline
    \textbf{Any Graph} & Sometimes & Sometimes & Sometimes & Sometimes \\
    \hline
    \begin{tabular}{@{}c@{}}{\bf Crossing Closed}\\ (Definition~\ref{crossingClosedDef})\end{tabular} & Sometimes & Always & Sometimes & Sometimes \\
    \hline
     \begin{tabular}{@{}c@{}}{\bf Upper Crossing Closed}\\ (Definition~\ref{upperCrossingClosedDef}) \end{tabular}& Sometimes & Always & Always & Sometimes \\
    \hline
    \begin{tabular}{@{}c@{}} {\bf Perfectly Labeled}\\ (Definition~\ref{perLabDef}) \end{tabular} & Always & Sometimes & Always & Always \\
    \hline
    \begin{tabular}{@{}c@{}} {\bf Tightly Closed}\\ (Definition~\ref{tightlyClosedDef}) \end{tabular} & Always & Always & Sometimes & Sometimes \\
    \hline
    \textbf{ \begin{tabular}{@{}c@{}}Upper Crossing Closed  \\ and Tightly Closed \end{tabular}} & Always & Always & Always & Always\\
    \hline
    \begin{tabular}{@{}c@{}}{ \bf Strongly Upper Crossed}\\ (Definition~\ref{stronglyUpperCrossedDef})\end{tabular}   & Always & Sometimes & Always & Always \\
       \hline

    \end{tabular}
    \caption{Families of graphs and their respective properties.}
  \label{MetaThmTable}
\end{table}

To finish this section, we gather all the information about families of graphs that we have seen throughout this paper.  This data appears in Table~\ref{MetaThmTable}. The rows of Table~\ref{MetaThmTable} refer to the families of graphs and the columns to the properties of the graphs or their noncrossing bond poset. The term ``NCNBC interpretation" refers to if the M\"obius function and characteristic polynomial (if applicable) have the noncrossing NBC set interpretation of Theorem \ref{NCNBCThm} or not. Every instance of ``sometimes" is genuine in the sense that there are graphs in that family which do and do not posses the prescribed property.

\section{Open Problems}\label{openProblemsSec}

As we have seen, several of the nice properties of the noncrossing partition lattice and the bond lattice have analogues in the noncrossing bond poset.  Given the multitude of nice properties that these lattices enjoy, we encourage the reader to see if their favorite properties have an analogue in the noncrossing bond poset.    We collect a few open problems that we have found interesting  below.  The list is in no way to be considered complete.

Recall that the Whitney numbers of the first kind of a graded poset are the numbers $w_0,w_1,\dots, w_n$ where $w_i$ is the sum of the M\"obius values of elements of $P$ of rank $i$.  In other words, they are the coefficients of the characteristic polynomial.  Moreover, recall that a sequence $a_0,a_1,\dots, a_n$ of real numbers is called \emph{log-concave} if for all $1\leq i\leq n-1$
we have that $a_{i-1}a_{i+1}\leq a_i^2$.

Gian-Carlo Rota conjectured that the Whitney numbers of the first kind for  geometric lattices (which include bond lattices) are log-concave.  In~\cite{h:mnphcp} Huh proved that the Whitney numbers of the first kind for bond lattices are log-concave and further work of Adiprasito, Huh, and Katz~\cite{ahk:htcg} proved the more general conjecture concerning the Whitney numbers of the first kind for geometric lattices.  Since the noncrossing bond poset is a (relatively) well-behaved subposet of a geometric lattice, it seems natural to ask if the corresponding conjectures hold for the noncrossing bond poset.  

\begin{question}
For which graphs are the Whitney numbers of the first kind of $\NC_G$ unimodal or log-concave?

\end{question}
\noindent  We should note that, unlike the case for the bond lattice, the Whitney numbers of the first kind of the noncrossing bond poset do not need to alternate in sign and can have internal zeros (e.g.~the characteristic polynomial of the 5-pointed star in Figure~\ref{5pointStarFig} has an internal zero).  As a result, it is not the case that the absolute values of the Whitney numbers of the first kind are log-concave or unimodal in general.  We mention this since, if the sequence did alternate and have no internal  zeros, the log-concavity would imply unimodality.  

 The noncrossing partition lattice is  well-known to be rank-symmetric (see, for example~\cite{k:nc}).  That is, for $\NC_{n+1}$, the number of elements of rank $k$ is the same as the number of elements of rank $n-k$.  It seems that it is rare for the noncrossing bond poset to be  rank-symmetric. This should not be that surprising as the bond lattice is also rarely rank-symmetric.  However, for $n\geq 5$, computations suggest that if we let $C_n$ denote the cycle on $n$ vertices with edges $12,23, \dots, n-1n, 1n$, then the complement $\overline{C_n}$  has a noncrossing bond poset which is rank-symmetric.   This leads us to a broader question.
\begin{question}
When is $\NC_G$ rank-symmetric?
\end{question}

As we saw in the  discussion preceding the previous question, the graph $\overline{C_n}$ seems to have a rank-symmetric noncrossing bond poset.  Despite this nice property, it seems that the poset is not shellable.    Naturally, this leads us to the following.

\begin{question}
For what graphs is the noncrossing bond poset shellable?  
\end{question}
\noindent We note here that $\overline{C_n}$ is tightly closed (but not upper crossing closed).  Thus, we know that tightly-closed (and hence crossing closed) does not imply shellability.  There is some hope that upper crossing closed graphs produce shellable noncrossing bond posets.  However, since they are not always graded, this will require considering non-pure shellings.

Given a graph (or more generally a matroid) one can consider the collection of non-broken circuits.  This set forms a simplicial complex called the \emph{broken circuit complex} or NBC complex.  It has several nice properties, its $f$-vector encode the coefficients of the chromatic polynomial of the graph and the complex is known to be shellable.  A related complex called the \emph{independence complex} is formed considering all the subsets of the edges sets which form acyclic subgraphs.  Since subsets of noncrossing sets are noncrossing, we can also consider the simplicial complex of noncrossing NBC sets and noncrossing independent sets.    
\begin{question}
What is the structure of the noncrossing NBC complex and noncrossing independence complex of a graph?
\end{question}

 \section*{Acknowledgements}
   Some of the results in this article had their genesis in the first author's master thesis.  We would like to thank Ed Allen, Hugh Howards, and Sarah Mason who made helpful suggestions and comments on that thesis. We would also like to thank the referee for their diligent reading and suggestions that greatly improved the article. We used {\tt SageMath}~\cite{sagemath} for computer explorations  and for generating conjectures.

\appendix
\section{Appendix: Algorithms}\label{algProofs}

In this appendix, we give two important algorithms on $\NC_G$, Algorithms \ref{crossingclosedalg} and \ref{uppercrossingclosedalg}. We also give a forbidden subgraph characterization  of upper crossing closed graphs in  Theorem \ref{ucccharacterization}.

Algorithm \ref{crossingclosedalg} decides if $\NC_G$ is a lattice in time on the order of $n^7$ where $n$ the number of vertices of the graph $G$. This is proved in Theorem~\ref{crossingclosedalgthm}.  Note that a brute-force algorithm to test if $NC_G$ is a lattice can take time super-exponential in $n$.  For example, an algorithm that checks if every pair of elements in $NC_G$ has a meet and a join, could take time at least on the order of the number of elements of $NC_G$.  For an $n$ vertex graph, that may be as large as the Bell number $B_n$ of the number of set partitions of $[n]$ and $B_n > (n/e \log(n))^n$ \cite{Berend}.  
Note that we call Algorithm \ref{crossingclosedalg} the ``crossing-closed'' algorithm as it is actually checking if $G$ is crossing closed.  Of course, $\NC_G$ being a lattice is equivalent to $G$ being crossing closed, see Theorem \ref{crossingClosedThm}. 

Algorithm \ref{uppercrossingclosedalg} determines if a graph $G$ is an upper crossing closed graph. Recall that if $G$ is upper crossing closed then the M\"obius function and characteristic polynomial of $\NC_G$ have nice interpretations in terms of noncrossing NBC sets, see Theorem \ref{NCNBCThm}. When given a graph $G$, Algorithm \ref{uppercrossingclosedalg} will either produce a specific upper crossing closed ordering of $E(G)$ or will produce what we term an obstruction (see Definition \ref{def:obstruction}), a specific subgraph of $G$ that clearly shows there can be no such ordering.  This also gives a forbidden subgraph characterization of upper crossing closed graphs, Theorem \ref{ucccharacterization}.

In Theorem \ref{uppercrossingclosedalgproof}, we prove that the Algorithm \ref{uppercrossingclosedalg} will run in time on the order of $n^8$ where again $n$ is the number of vertices of $G$. Note that a brute force algorithm could again take time super-exponential in $n$, if it is forced to test some positive fraction of the $\binom{n}{2}$! possible orderings on the edges of $G$.

We first present our algorithm that decides if $G$ is crossing closed, i.e.~if $\NC_G$ is a lattice.

\begin{algorithm} {\bf Crossing Closed Algorithm} \label{crossingclosedalg}

\noindent {\bf Input}: A graph $G$ on $[n]$.

\noindent {\bf Output}: A yes/no decision as to whether $G$ is crossing closed, or, equivalently, whether $\NC_G$ is a lattice.

\noindent {\bf Method:} 
For each pair of crossing edges $e$ and $f$ find a shortest path $P(e,f) = x_0 x_1 \dots x_k$
with $e = x_0x_1$ and $f = x_{k-1}x_k$
and $k \geq 3$. If for some crossing pair $e$ and $f$, $P(e,f)$ fails to exist or has $k \geq 4$ and has some vertex $x_i$ with $2 \leq i \leq k-2$ such that $x_i$ does not separate $e$ and $f$ then return ``No, $G$ is not crossing closed.'' Otherwise return ``Yes, $G$ is crossing closed.'' 
\end{algorithm}

\begin{thm} \label{crossingclosedalgthm}
Algorithm \ref{crossingclosedalg} is a correct algorithm that runs in time $O(n^7)$ where $n$ is the number of vertices of $G$.
\end{thm}

\begin{proof}
First, we will compute the complexity of the algorithm. The Floyd-Warshall algorithm gives a shortest path between all pairs of vertices in $O(n^3)$ time \cite{Floyd}.  With that pre-processing done, there are at most $\displaystyle\binom{n}{4}$  pairs of crossing edges to check.  For each pair $e$ and $f$ of crossing edges, there are at most $n$ vertices on the shortest path connecting them to check.  Checking that one of those vertices separates $e$ and $f$ can be done by breadth-first search in $O(n^2)$ time so this algorithm will run in $O(n^7)$ time.

Next, we show that the algorithm always gives the correct output. Suppose $G$ is crossing closed.  We will show that the algorithm will return a ``yes''.  For every pair of crossing edges $e$ and $f$, $J(e,f)$ exists.  If there is an edge incident to $e$ and $f$ (where incident means having a common vertex), then $P(e,f)$ will be a path $(e,g,f)$ for some edge $g$ connecting $e$ and $f$ and the algorithm will not give a ``no'' answer based on this pair.  If there is no edge incident to both $e$ and $f$, then by Lemma \ref{mefLemma}, $J(e,f)$ will be a dumbbell graph, see Figure \ref{dumbellGraphFig}, a graph induced by $e, f$ and a path $Q$ with one end adjacent to $e$ and the other end adjacent to $f$ with all vertices in $Q$ separating $e$ and $f$.  Thus the vertices $\{x_2, \dots, x_{k-2}\}$ on $P(e,f)$ must be the vertices of $Q$. The algorithm will not give a ``no'' answer based on this pair either. Thus the algorithm will return a ``yes''.

Suppose now that the algorithm returns a ``yes''.  We will show that $G$ is crossing closed.  Let $e$ and $f$ be a pair of crossing edges.  Since the algorithm returned a ``yes'', a shortest path $P(e,f)$ must exist. If the path has 3 vertices, i.e.~there is an edge connecting $e$ and $f$, then $J(e,f) = G[e \cup f]$ exists.  Suppose now that the path contains at least 4 vertices. Let $M$ be the subgraph of $G$ induced on the vertices of $P(e,f)$.  We claim that $M$ is contained in every connected, induced subgraph of $G$ that contains $e$ and $f$ and so $J(e,f)=M$.  Let $x$ be any vertex in $P(e,f)$ that is not in $e$ or $f$ and let $H$ be an connected induced subgraph of $G$ containing $e$ and $f$.  If $x$ is not in $H$ then it cannot separate $e$ and $f$.  Since the algorithm returned ``yes'', $x$ must separate $e$ and $f$ and so $x$ must be in $H$.  Thus every vertex of $M$ is in $H$ and since they are both induced, $M$ is in $H$.
\end{proof}

Now we turn our attention to the problem of deciding whether a graph is upper crossing closed. First, let us note that not all crossing closed graphs are upper crossing closed.  As an example, consider the 5-pointed star in Figure~\ref{5pointStarFig}.  It is not hard to verify that if $e$ and $f$ cross in the 5-pointed star, $J(e,f)$ is a subgraph of $K_4$ and so is crossing closed.  However, every edge of the graph is crossed and so it is impossible to have an ordering that is upper crossing closed as the smallest edge must be noncrossing.  As it turns out, this kind of issue is the only obstacle to a crossing closed graph being upper crossing closed.

\begin{definition} \label{def:obstruction}
Let $G$ be a crossing closed graph.  We say a subgraph $H$ of $G$ is an  \emph{obstruction} to $G$ being upper crossing closed if for  every edge $e$ in $H$ there is an edge $f$ in $H$ which crosses $e$ such that  $J(e,f) \subseteq H$.  
\end{definition}

Theorem \ref{uppercrossingclosedalgproof} proves that Algorithm \ref{uppercrossingclosedalg} below will, when given a graph $G$, either produce an upper crossing closed ordering on $E(G)$ or an obstruction.  It also proves that an obstruction demonstrates no such ordering is possible.  Thus  we get the following structural characterization of upper crossing closed graphs.

\begin{thm} \label{ucccharacterization}
A graph $G$ is upper crossing closed if and only if it contains no obstruction as a subgraph.
\end{thm}

Note also that if every edge $e$ in $G$ crosses some other edge of $G$, then $G$ itself is an obstruction of $G$.  Thus we also have the following.
\begin{cor} \label{ucccharacterizationcor}
If $G$ is a graph with every edge crossing some other edge, then $G$ is not upper crossing closed.
\end{cor}
Note that Corollary \ref{ucccharacterizationcor} shows that the graph $G$ of Figure \ref{nonGradedUpperCCFig} and the 5-pointed star of Figure \ref{5pointStarFig} are not upper crossing closed.

\begin{algorithm} {\bf Upper Crossing Closed Algorithm} \label{uppercrossingclosedalg}

\noindent {\bf Input:} A graph $G$ on $[n]$.

\noindent {\bf Output:} A yes/no decision on whether $G$ is crossing closed. Then if $G$ is crossing closed, a yes/no decision on whether $G$ is upper crossing closed.  If $G$ is upper crossing closed, an upper crossing closed ordering is produced, and if $G$ is crossing closed but not upper crossing closed, an obstruction is produced.

\noindent {\bf Method:}
\begin{enumerate}
\item Run the crossing closed algorithm on $G$, Algorithm \ref{crossingclosedalg}. If the answer is no, return ``No.  $G$ is not crossing closed and hence not upper crossing closed." and terminate.  If the answer is yes, return ``Yes.  $G$ is crossing closed." and continue.
    
\item Let $L = \emptyset$ and let $\sigma = \emptyset$.  (Throughout the algorithm, $L$ will be a subset of $E(G)$ and $\sigma$ will be an ordering on $L$.)

\item Let $L'$ be the set of edges $e$ in $E(G) \setminus L$ such that for every edge $f \in E(G) \setminus L$ that crosses $e$, $E(J(e,f)) \cap L \neq \emptyset$.

\item If $L' \neq \emptyset$ update $L$ to be $L \cup L'$ and update $\sigma$ to be the ordering
on $L \cup L'$ that orders $L$ according to $\sigma$
and then puts all the edges of $L'$ after the edges of $L$.  The ordering within $L'$ can be arbitrary.  Go back to step 3.

\item If $L' = \emptyset$, decide on the output of the algorithm. If $L = E(G)$, return ``Yes, $G$ is upper crossing closed, and $\sigma$ is an upper crossing closed ordering on $E(G)$.''  If $L \neq E(G)$, return ``No. $G$ is not upper crossing closed, and the spanning subgraph of $G$ with edge set $E(G) \setminus L$ is an obstruction.''

\end{enumerate}

\end{algorithm}

We will now show how the algorithm runs on two graphs, one upper crossing closed and the other not. First, let $G$ be the twisted 4-cycle in Figure~\ref{twistedC4Fig}. As we have already seen, $G$ is crossing closed and upper crossing closed with respect to the lexicographic order on its edges. The algorithm will thus correctly conclude that $G$ is crossing closed in step 1 and will set $L=\emptyset$ and $\sigma=\emptyset$ in step 2. Next, it will go to step 3.  Since $L$ is empty, $L'$ is the set of edges which cross no other edges.  So $L'=\{12,34\}$.  Then the algorithm passes to step 4 where $L$ is set to be $\{12,34\}$ and $\sigma$ is set to be some total ordering of $\{12,34\}$.  Now we return to step 3.  Now, $E(G) \setminus L = \{13, 24\}$. Since $13$ and $24$ form the only crossing in $G$ and $J(13,24) = G$ intersects $L$, the algorithm sets $L' = \{13,24\}$. Next, we go to step 4, where $L$ is set to be $\{12, 34, 13,24\}$ and $\sigma$ is some total ordering where the first two elements are $12$ and $34$ and the last two elements are $13$ and $24$.   Then we return to step 3, where $L'$ is set to be empty.  Finally, we go to step 5 and since $L=E(G)$, the algorithm returns that $G$ is upper crossing closed with respect to the ordering $\sigma$, which indeed it is. The reader may have noticed that the ordering the algorithm produces is not the lexicographic ordering. This is because the algorithm always puts edges with no crossing before any edge with a crossing.  Thus the algorithm is not always capable of producing all possible upper crossing closed orderings.

Now we give an example of how the algorithm runs on a graph $G$ that is not upper crossing closed. Let $G$ be the 5-pointed star in Figure~\ref{5pointStarFig}. Since $G$ is crossing closed, the algorithm will pass to step 2 and set $L=\emptyset$ and $\sigma =\emptyset$.  Then it moves to step 3.  Since all the edges in $G$ cross some other edge and $L$ is empty, $L'$ is empty too. As a result, the algorithm moves to step 5.  Since $L\neq E(G)$, the algorithm returns that $G$ is not upper crossing closed and correctly provides $G$ as an obstruction. We should note that for any graph in which every edge crosses another edge, the algorithm will terminate with $L'= L = \emptyset$ and thus will correctly conclude that the spanning subgraph with edge set $E(G) \setminus L = E(G)$, i.e.~$G$ itself, is an obstruction.

We now prove that Algorithm \ref{uppercrossingclosedalg} is correct and runs in polynomial time.

 \begin{thm} \label{uppercrossingclosedalgproof}
Let $G$ be a graph.  Then we have the following.
\begin{enumerate}
    \item[(a)] Algorithm \ref{uppercrossingclosedalg} runs in time $O(n^8)$ where $n$ is the number of vertices of $G$.
   
    \item[(b)] 
    If Algorithm \ref{uppercrossingclosedalg} concludes by giving a purported obstruction $H$, then $H$ is indeed an obstruction.
    \item[(c)]  If Algorithm \ref{uppercrossingclosedalg} produces an obstruction, then G is not upper crossing closed.
    \item[(d)]   If Algorithm \ref{uppercrossingclosedalg} does not produce an obstruction, then $G$ is upper crossing closed and the order $\sigma$ it produces is an upper crossing closed ordering.
    \item[(e)]    Algorithm \ref{uppercrossingclosedalg} is correct.
\end{enumerate}

\end{thm} 
 \begin{proof}
 First, we show (a).  The crossing closed algorithm runs in $O(n^7)$ time as a subroutine in step 1. During the course of that run, shortest paths $P(e,f)$ connecting all pairs of crossing edges $e$ and $f$ are created.  By the proof of Theorem \ref{crossingclosedalgthm}, these paths determine $J(e,f)$ for each such pair of edges.  When running step 3, there are at most $n^2$ edges $e \in E(G) \setminus L$ to check and then for each such $e$ there are at most $n^2$ edges $f$ that cross $e$ to check. Since the $J(e,f)$ are already calculated it takes $n^2$ comparisons to calculate $J(e,f) \cap L$, so step 3 takes $O(n^6)$ time each time it is run.  It is run at most $n^2$ times so the algorithm takes $O(n^8)$ time.
 
 Now we verify (b).  If the algorithm terminates with $L' = \emptyset$ and $L \neq E(G)$, then the output of the algorithm is the spanning subgraph $H$ with edge set $E(G) \setminus L$.  Since $L' = \emptyset$, every edge $e \in E(G) \setminus L$ must cross another edge $f \in E(G) \setminus L$ such that $E(J(e,f)) \subseteq E(G) \setminus L$.  Thus $H$ is indeed an obstruction.
 
 Now we verify (c). Let $H$ be an obstruction. Then for every ordering $\unlhd$ of $E(G)$, the first edge $e$ of $H$ will cross some other edge $f$ of $H$ with $J(e,f) \subseteq H$.  But since $J(e,f) \subseteq H$ and $e$ is the minimum edge of $H$, no edge $g \in J(e,f)$ will satisfy $g \lhd e,f$.  It follows that $G$ is not upper crossing closed with respect to any ordering.
 
 Next, we prove (d). Suppose no obstruction is found. We claim that the ordering $\sigma$ on $E(G)$ that is produced is an upper crossing closed ordering.  Let $e$ and $f$ be a pair of crossing edges in $G$.  Consider the first point in time during the run of the algorithm in which $e,f \not \in L$ and $e$ or $f$ or both are in $L'$.  Say $e \in L'$.  Then $J(e,f) \cap L \neq \emptyset$ otherwise $e$ would not be in $L'$.  So there will be an edge $g \in J(e,f) \cap L$.  In the ordering $\sigma$, all the edges in $L$ are less than all the edges not in $L$ so $g$ will be less than $e$ and $f$. This shows that $G$ is upper crossing closed with respect to $\sigma$.
 
 Finally, let us show (e).  Suppose that $G$ is not upper crossing closed.  Then the algorithm must find an obstruction. If it did not, by part (d), $\sigma$ would be an upper crossing closed ordering. By part (b), what the algorithm produces is really an obstruction and by part (c) this obstruction demonstrates that $G$ is not upper crossing closed. Thus the algorithm will return an obstruction and correctly returns that $G$ is not upper crossing closed.
 
 Now suppose that $G$ is upper crossing closed.  Then by the contrapositive of part (c) and by part (b) the algorithm produces no purported obstruction.  So then by part (d), it returns an upper crossing closed ordering.  It will then correctly return that the graph is upper crossing closed.
 \end{proof}

\end{document}